\documentclass[sn-mathphys,Numbered]{sn-jnl}

\usepackage{graphicx}%
\usepackage{multirow}%
\usepackage{amsmath,amssymb,amsfonts, scalefnt}%
\usepackage{amsthm}%
\usepackage{mathrsfs}%
\usepackage[title]{appendix}%
\usepackage{xcolor}%
\usepackage{textcomp}%
\usepackage{manyfoot}%
\usepackage{booktabs}%
\usepackage{algorithm}%
\usepackage{algorithmicx}%
\usepackage{algpseudocode}%
\usepackage{listings}%
\usepackage{natbib}
\usepackage[export]{adjustbox}

\algnewcommand\algorithmicforeach{\textbf{for each}}
\algdef{S}[FOR]{ForEach}[1]{\algorithmicforeach\ #1\ \algorithmicdo}

\theoremstyle{thmstyleone}%
\theoremstyle{thmstyletwo}%

\theoremstyle{thmstylethree}%

\begin{document}

\title[Article Title]{Cislunar Satellite Constellation Design Via Integer Linear Programming\footnote{This paper is a substantially revised version of Paper AAS 23-358 presented at the AAS/AIAA Astrodynamics Specialist Conference, Big Sky, MT, in Aug. 2023. The added content includes new analysis, case studies, and results.}}


\author[1]{\fnm{Malav} \sur{Patel}}\email{mpatel636@gatech.edu}

\author[1]{\fnm{Yuri} \sur{Shimane}}\email{yshimane3@gatech.edu}

\author[2]{\fnm{Hang Woon} \sur{Lee}}\email{hangwoon.lee@mail.wvu.edu}

\author*[1]{\fnm{Koki} \sur{Ho}}\email{kokiho@gatech.edu}

\affil*[1]{\orgdiv{Aerospace Engineering}, \orgname{Georgia Institute of Technology}, \orgaddress{\street{270 Ferst Drive}, \city{Atlanta}, \postcode{30332}, \state{Georgia}, \country{USA}}}

\affil[2]{\orgdiv{Mechanical and Aerospace Engineering}, \orgname{West Virginia University}, \orgaddress{\street{1306 Evansdale Drive}, \city{Morgantown}, \postcode{26506}, \state{West Virginia}, \country{USA}}}

\abstract{Cislunar space domain awareness is of increasing interest to the international community as Earth-Moon traffic is projected to increase, which raises the problem of placing space-based sensors optimally in a constellation to satisfy the space domain awareness demand in cislunar space. This demand profile can vary over space and time, making the design optimization problem challenging. This paper tackles the problem of satellite constellation design for spatio-temporally varying coverage demand by leveraging an integer linear programming formulation. The developed optimization formulation assumes the circular restricted 3-body dynamics and attempts to minimize the number of satellites required for the requested demand profile.}


\maketitle

\section{Introduction}\label{sec1}

The Earth-Moon system and related orbital trajectories have been of interest since the early Apollo missions in the 1960s and 1970s. This trend carried into the early part of the 21st century with research efforts aiming to find better trajectories for tours in the Earth-Moon system and operations near the Moon \cite{folta2013preliminary, haapala2013trajectory}.

The recent lunar exploration activities by government and private sectors have renewed this interest as Earth-Moon traffic is projected to increase through crewed and uncrewed missions in the coming decade. This makes the region of space between the Earth and the Moon, coined \textit{Cislunar Space}, an attractive region to monitor. Ground-based optical measurement is unable to provide sufficient observational power for such a vast region, particularly for the regions near the direction of the Moon. This insufficiency motivates space-based monitoring schemes, realized as satellite constellations with optical payloads on favorable orbits for surveillance. 

Following this demand for surveillance capability, various studies examined the orbit selection problems for cislunar SDA \cite{cunio2020payload, thompson2021cislunar}. Wilmer et al. \cite{wilmer} and Dahlke et al. \cite{dahlke2022preliminary} investigated using periodic orbits to monitor activity near the Lagrange points (L1-L5) in the Earth-Moon system. Zimovan et al. \cite{zimovan2017near} investigated near-rectilinear halo orbits (NRHOs) around the Moon because of attractive stability, transfer prospectives, and eclipsing properties. Because of this, NRHOs are promising orbits for lunar space stations that serve as outposts for Cislunar activity. Frueh et al. \cite{frueh2021cislunar} introduced the 2:1 resonant orbit as an attractive candidate for monitoring cislunar space because of its close approaches to the Earth and Moon and because it covers the entire cislunar region in under 20 revolutions. Gupta et al. \cite{guptaconstellation} investigated how phasing an observer and a data relay satellite on separate orbits can impact the length of the window in which communication is possible. 

With those orbit options, studies have also focused on optimizing the space-based cislunar SDA architecture. Klonowski et al. \cite{Klonowski} formulated SDA architecture design as a Markov Decision Process and employed the popular Monte Carlo Tree Search to optimize a constellation of up to 5 observers to maintain custody of a space object traveling on a specific orbit between the Earth and Moon. Vendl et al. \cite{vendl2021cislunar} explored a more general \textit{region} of cislunar space, specifically a 20-degree cone emanating from the Earth towards the Moon. The region consisted of 1071 target points, which were used to choose an optimal phase of an observer on an L1 (and separately, L2) Lyapunov orbit to maximize observation capability. However, the number of satellites was not a design variable in this formulation. In contrast, Visonneau et al. \cite{visonneau2023optimizing} presented a genetic algorithm for finding optimally sized constellations for monitoring regions of cislunar space, with no restriction on number of observers.

Although these aforementioned works tackled the space-based cislunar SDA architecture problems, they do not provide a method to optimize the cislunar SDA architecture to satisfy the spatio-temporally varying demand effectively. This need is critical for cislunar SDA because we often have a set of times and regions of particular interest rather than constantly covering the entire cislunar space. 

This paper tackles this challenge by proposing a new integer linear programming (ILP) formulation. This formulation is inspired by the previous work by Lee et al. \cite{2020JSpRo..57.1309L}, which developed the Access-profile-Pattern-Coverage (APC) decomposition method that formulates an ILP problem to optimize an Earth-based constellation that satisfies a temporally varying (e.g., day-time demand higher than night-time demand) and spatially varying demand (e.g., certain regions of higher interest than other regions).  By "demand", we mean that a target point is tasked with being observed at a specified time by a specified number of observer satellites. This work extends this formulation to cislunar SDA applications using the circular restricted three-body problem (CR3BP) dynamics to identify the optimal satellite orbit (including the phasing) to satisfy the given demand with the minimum number of satellites. As an example, the demand is represented as one or more moving targets in cislunar space, but the method is general enough to handle any spatio-temporally varying demand.

The rest of the paper is structured in the following way: Section \ref{sec2} introduces the generic ILP problem formulation based on APC decomposition, Section \ref{sec3} extends the formulation into cislunar space, and Section \ref{sec4} shows the analysis results, and Section \ref{sec5} concludes the paper.

\section{APC Decomposition and ILP Problem Formulation}\label{sec2}

The problem of interest to this paper is to minimize the number of satellites required to satisfy a given spatio-temporally varying coverage requirement. We will introduce the problem formulation as an ILP problem based on the APC decomposition method developed by Lee et al. \cite{2020JSpRo..57.1309L}. This section introduces the generic ILP formulation resulting from the APC decomposition method, which will be later extended to the cislunar SDA problem in Section \ref{sec3}.
 
\subsection{Accessibility Profile}\label{subsec2}
The accessibility profile is a vector of binary variables describing when the target is in view of the observer as the observer travels along its orbit. 

Suppose there is a target point's position vector in a given reference frame, denoted by $\boldsymbol{p}$. Suppose we have an observer at position $\boldsymbol{r}$. The target point's observability can be quantified by an accessibility measure, $g$. A simple measure could be the distance between observer and target, $g(\boldsymbol{p}, \boldsymbol{r}) = \|\boldsymbol{p} - \boldsymbol{r}\|$. A more sophisticated measure is the apparent magnitude of the object, which is a function of the position of the sun as well, $g = g(\boldsymbol{p}, \boldsymbol{r}, \boldsymbol{s})$ 
\cite{1974STIN...7512024K}.

Suppose the observer is on an orbit with period $T$, time is discretized for computational ease, and $\boldsymbol{p}$ is a static target point in space. 
At time step $n$, $\boldsymbol{r}_n$ is the observer's position, $\boldsymbol{s}_n$ the sun's position,  and $a[n] = g(\boldsymbol{p}, \boldsymbol{r}_n, \boldsymbol{s}_n)$ is the accessibility measure. If the orbit is discretized into $L$ equal time steps, then $\boldsymbol{a}$ is a $L \times 1$ vector containing the accessibility measures along the observer's orbit. We assume that the target point will be in one of two states, [observable] or [not observable]. In other words, the continuous accessibility measure is replaced by a binary one, \textit{accessibility profile}, $\boldsymbol{v}$. $\boldsymbol{a}$ is converted into $\boldsymbol{v}$ by using a threshold value $M$ (e.g., threshold for the apparent magnitude of the target):

\begin{equation}
v[n] = 
    \begin{cases}
        1, & \text{if } a[n] \leq M\\
        0, & \text{otherwise}
    \end{cases}
\end{equation}
At time step $n$, $v[n]$ indicates whether the observer can view the target point.

\subsection{Constellation Pattern Vector}\label{subsec2}
Instead of a single observer, suppose we have a set (called a constellation) of observers on the same orbit of period $T$ again discretized into L time steps. The \textit{constellation pattern vector} $\boldsymbol{x}$ captures the phasing of the observers in the constellation relative to a \textit{seed satellite}. Suppose that at time step $n=1$, the seed satellite is at position $\boldsymbol{r}$, and after $n_k$ time steps, the $k$-th observer reaches position $\boldsymbol{r}$. In this case, the time delay for the $k$-th observer is denoted as $\Delta t_k = (n_k-1) \times \frac{T}{L}$. The constellation pattern is defined as a vector $\boldsymbol{x}$, where $x[n]=1$ if $n=n_k$ for some $k$ and $0$ otherwise. 

\subsection{Coverage Timeline and Coverage Requirement}\label{subsec2}
The \textit{coverage timeline}, denoted by $\boldsymbol{b}$, is a $L \times 1$ vector that extends the accessibility profile concept to multiple observers. At time step $n$, $b[n]$ describes the number of satellites with view access to the target point. Instead of a vector of binary variables, the coverage timeline is a vector of non-negative integers.

 Interestingly, Lee et al. \cite{2020JSpRo..57.1309L} show that a circular convolution between the constellation pattern vector $\boldsymbol{x}$ and the accessibility profile $\boldsymbol{v}$ returns the coverage timeline $\boldsymbol{b}$. The circular nature of the convolution is a direct consequence of the periodicity of the observers' orbit. Note that periodicity is a requirement for the circular convolution operation to be valid. Also, note that solar illumination conditions may vary after each period if the orbits are not perfectly resonant, in which case $\boldsymbol{b}$ would vary after each period. We acknowledge this in Sec. \ref{subsec4} and show that differences in solar illumination conditions could still lead to satisfactory performance.  

The circular convolution can be conveniently expressed as a matrix multiplication between an $L \times L $ circulant matrix, $V$, and the constellation pattern vector, $\boldsymbol{x}$. Here, the $i$-th column of $V$ is produced by circularly shifting $v$ by $i$ elements.

\begin{equation}
    V\boldsymbol{x} = \boldsymbol{b}
\end{equation}

Expanded in matrix form and using 1-based indexing,

\begin{equation}
\scalefont{1.0}
{
\begin{aligned}
{\begin{bmatrix}v[1]&v[L]&\cdots &v[3]&v[2]\\v[2]&v[1]&v[L]&\hdots&v[3]\\v[3] &v[2]&v[1]& &\vdots \\\vdots&\vdots&\ddots &\ddots &\\v[L]&v[L-1]&\cdots &&v[1]\\\end{bmatrix}}
{\begin{bmatrix}
    x[1]\\x[2]\\\vdots\\\\x[{L}]
\end{bmatrix}}
= 
{\begin{bmatrix}
    b[1]\\b[2]\\\vdots\\\\\b[{L}]
\end{bmatrix}}
\end{aligned} 
}
\end{equation}

Often, we are interested in satisfying a coverage requirement, $\boldsymbol{f}$, which demands a different number of observers at different time steps. The goal, then, is to place a set of observers on the orbit (i.e. find a constellation pattern $\boldsymbol{x}$) such that the resulting coverage timeline $\boldsymbol{b}$ meets or exceeds the coverage requirement $\boldsymbol{f}$.
Stated mathematically,

\begin{equation}
\begin{aligned}
& {\text{Find}}
& & \boldsymbol{x} \\
& \text{such that}
& & b[n] \geq f[n] &&& \text{for}\  n = \{1, 2, \cdots, L\}\\
& \text{where} 
&& V\boldsymbol{x} = \boldsymbol{b}
\end{aligned}
\end{equation}

\subsection{ILP Formulation}\label{subsec2}
In general, there can be many constellation pattern vectors that satisfy the coverage requirement. However, it is often beneficial to minimize the total number of observers (or another metric) in the constellation while still satisfying the coverage requirement, since fewer observers reduce launch and maintenance costs.

The problem above of finding a constellation pattern $\boldsymbol{x}$ is now a problem of finding the \textit{optimal} constellation pattern $\boldsymbol{x}^*$ that minimizes the number of observers in the constellation. The following integer linear program presents this mathematically.

\begin{equation}
\begin{aligned}
& \underset{\boldsymbol{x}}{\text{minimize}}
& & \boldsymbol{1}^T \boldsymbol{x} \\
& \text{subject to}
& & V\boldsymbol{x} \geq \boldsymbol{f} \\
&&& \boldsymbol{x} \in \{0 ,1\}^L
\end{aligned}
\end{equation}

 The objective is to minimize the total number of satellites, hence the objective $\boldsymbol{1}^T \boldsymbol{x}$, which sums all the entries in $\boldsymbol{x}$. 

\subsection{Extension to Multiple Orbits \& Target Points}\label{subsec2}
The above formulation is linear, and thus allows extensions to multiple target points and multiple constellations.
\subsubsection{Multiple Target Points}\label{subsubsec1}
Suppose we have a set $\mathcal{J}$ of target points and we seek to find an optimal constellation of observers on a single orbit to monitor these target points according to a set of coverage requirements $\mathcal{F} = \{\boldsymbol{f}_1, \boldsymbol{f}_2, \cdots, \boldsymbol{f}_{\lvert\mathcal{J}\rvert}\}$, where $\lvert\mathcal{J}\rvert$ is the number of target points. Each target point in $\mathcal{J}$ will have its own accessibility profile, and thus its own circulant matrix. Denote the circulant matrix for each target point by $V_{j,1}$, where $j$ indexes the elements in $\mathcal{J}$.  Then we have

\begin{equation}
\begin{aligned}
& \underset{x}{\text{minimize}}
& & \boldsymbol{1}^T x \\
& \text{subject to}
& & V_{1,1}\ \boldsymbol{x} \geq \boldsymbol{f}_1 \\
& & & V_{2,1}\ \boldsymbol{x} \geq \boldsymbol{f}_2 \\
& & & \vdots \\
& & & V_{\lvert\mathcal{J}\rvert,1}\ \boldsymbol{x} \geq \boldsymbol{f}_{\lvert\mathcal{J}\rvert} \\
&&& \boldsymbol{x} \in \{0 ,1\}^L
\end{aligned}
\end{equation}

The linear constraint can be condensed into a single matrix equation by vertically concatenating the circulant matrices $V_{j,1}$ and coverage requirements $\boldsymbol{f}_j$ as follows:

\begin{equation}
\begin{aligned}
& \underset{x}{\text{minimize}}
& & \boldsymbol{1}^T \boldsymbol{x} \\
& \text{subject to}
&& {\begin{bmatrix}
    V_{1,1}\\V_{2,1}\\\vdots\\V_{\lvert\mathcal{J}\rvert,1}
\end{bmatrix}}
\boldsymbol{x}
\geq
{\begin{bmatrix}
    \boldsymbol{f}_1\\\boldsymbol{f}_2\\\vdots\\\boldsymbol{f}_{\lvert\mathcal{J}\rvert}
\end{bmatrix}} \\ 
&&& \boldsymbol{x} \in \{0 ,1\}^L
\end{aligned}
\end{equation}

\subsubsection{Multiple Constellation Orbits}\label{2}
Suppose we have a set $\mathcal{Z}$ of orbits and we seek to an optimal constellation of observers on each orbit to monitor a single target point. Each orbit in $\mathcal{Z}$ will have its own accessibility profile, and thus its own circulant matrix. Denote the circulant matrix for each orbit by $V_{1,z}$ where $z$ indexes the elements in $\mathcal{Z}$. Furthermore, each orbit will have its own constellation pattern vector, $\boldsymbol{x}_z$. Then we have

\begin{equation}
\begin{aligned}
& \underset{x}{\text{minimize}}
& & \sum_{z=1}^{\lvert\mathcal{Z}\rvert}  \boldsymbol{1}^T \boldsymbol{x}_z \\
& \text{subject to}
& &  V_{1,1}\ \boldsymbol{x}_1 \geq \boldsymbol{f}\\& & &  V_{1,2}\ \boldsymbol{x}_2 \geq \boldsymbol{f}\\&&& \ \vdots\\ &&&  V_{1,\lvert\mathcal{Z}\rvert}\ \boldsymbol{x}_{\lvert\mathcal{Z}\rvert} \geq \boldsymbol{f} \\
&&& x \in \{0 ,1\}^L
\end{aligned}
\end{equation}
where $\lvert\mathcal{Z}\rvert$ is the number of candidate orbits

The linear constraint can be condensed into a single matrix equation by horizontally concatenating the circulant matrices $V_{1,z}$ and vertically concatenating the constellation pattern vectors $\boldsymbol{x}_z$.

\begin{equation}
\begin{aligned}
& \underset{\boldsymbol{x}}{\text{minimize}}
& & \sum_{z=1}^{\lvert\mathcal{Z}\rvert}  \boldsymbol{1}^T \boldsymbol{x}_z \\
& \text{subject to}
&& {\begin{bmatrix}
    V_{1,1}&V_{1,2}&\cdots&V_{1,\lvert\mathcal{Z}\rvert}
\end{bmatrix}}
{\begin{bmatrix}
    \boldsymbol{x}_1\\\boldsymbol{x}_2\\\vdots\\\boldsymbol{x}_{\lvert\mathcal{Z}\rvert}
\end{bmatrix}}
\geq \boldsymbol{f}\\ 
&&& \boldsymbol{x}_z \in \{0 ,1\}^L \quad \forall z
\end{aligned}
\end{equation}

We can combine the linearities above to provide a single comprehensive ILP. Let $\boldsymbol{x}$ be the augmented constellation pattern vector obtained by vertically concatenating the constellation pattern vectors $x_z$, where $z$ indexes the elements in $\mathcal{Z}$. Let $\boldsymbol{f}$ be the augmented coverage requirement obtained by vertically concatenating the coverage requirements $\boldsymbol{f}_j$ for each target point, where $j$ indexes the elements in $\mathcal{J}$. Let $\boldsymbol{V}$ be the augmented circulant matrix obtained by vertically concatenating the circulant matrices for each target point and horizontally concatenating the circulant matrices for each orbit. For clarity, $\boldsymbol{V}$, $\boldsymbol{x}$, and $\boldsymbol{f}$ are shown in matrix form,
\begin{equation}
\begin{aligned}
\boldsymbol{V}
=
{\begin{bmatrix}V_{1, 1}&V_{1,2}&\cdots & \cdots &V_{1,\lvert\mathcal{Z}\rvert}\\V_{2, 1}&V_{2,2}&\cdots&\cdots&V_{2, \lvert\mathcal{Z}\rvert}\\\vdots &\vdots&\ddots&\ddots &\vdots \\\vdots& \vdots&\ddots &\ddots &\vdots\\V_{\lvert\mathcal{J}\rvert, 1}&V_{\lvert\mathcal{J}\rvert, 2}&\cdots &\cdots&V_{\lvert\mathcal{J}\rvert, \lvert\mathcal{Z}\rvert}\\\end{bmatrix}}, 
\quad
\quad
 \boldsymbol{x}
=
{\begin{bmatrix}
    \boldsymbol{x}_1\\\boldsymbol{x}_2\\\vdots\\\boldsymbol{x}_{\lvert\mathcal{Z}\rvert}
\end{bmatrix}}, 
\quad
\quad
\boldsymbol{f}
=
{\begin{bmatrix}
    \boldsymbol{f}_1\\\boldsymbol{f}_2\\\vdots\\\boldsymbol{f}_{\lvert\mathcal{J}\rvert}
\end{bmatrix}}
\end{aligned}
\end{equation}
Thus, our full ILP formulation is written as

\begin{equation}
\begin{aligned}
& \underset{x}{\text{minimize}}
& & \boldsymbol{1}^T \boldsymbol{x} \\
& \text{subject to}
& & \boldsymbol{Vx} \geq \boldsymbol{f}, \\
&&& \boldsymbol{x} \in \{0 ,1\}^{\lvert\mathcal{Z}\rvert L}
\end{aligned}
\end{equation}
where $\boldsymbol{1}^T$ has dimensions conformable to those of $\boldsymbol{x}$. This optimization method can be solved by commercial solvers like Gurobi.

\section{ILP Formulation for Cislunar SDA Architecture Design}\label{sec3}
\begin{figure}[H]
	\centering\includegraphics[width=\linewidth]{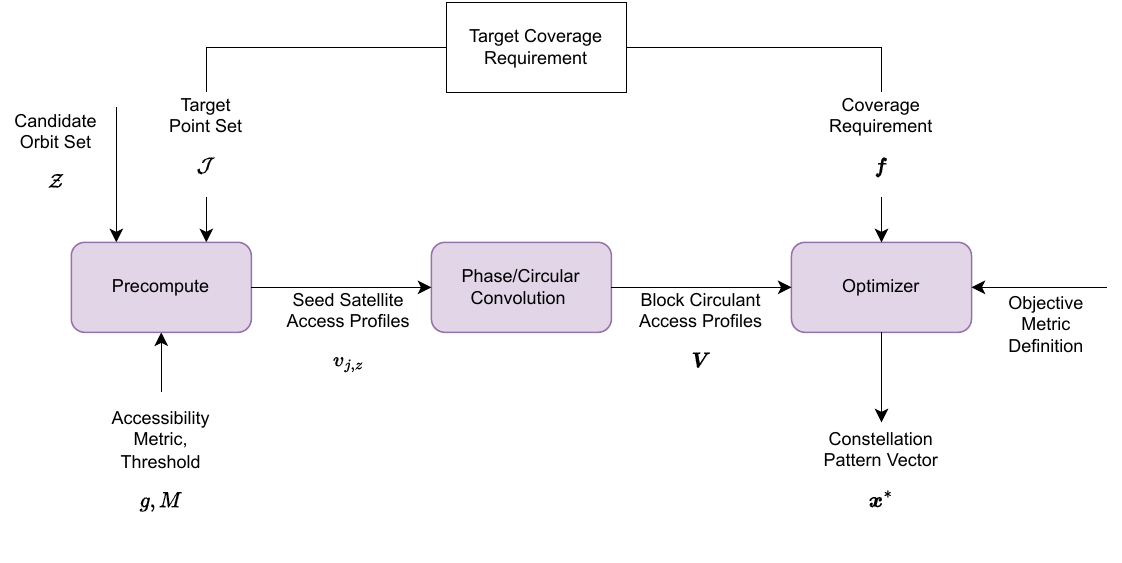}
	\caption{Optimizer Illustration. Purple boxes correspond to computation pipeline units. Arrows pointing into these boxes and associated text are inputs to said boxes. Arrows pointing out of these boxes and associated text are outputs of said boxes. In the white box is the Target Coverage Requirement which is used to produce two inputs into separate units of the computation pipeline.}
	\label{fig:2}
\end{figure}

This section extends the generic formulation introduced in Section \ref{sec2} to the cislunar SDA architecture design problem. The flow graph of the proposed formulation is shown in Figure \ref{fig:2}. The first stage consists of pre-computation, where we determine the accessibility for each \{target point, observer location\} pair in the entire constellation for the simulation period. Second, we construct the augmented circulant matrix $\boldsymbol{V}$ by phasing the accessibility measures. This primarily consists of a simple circular shift operation. Lastly, we provide the necessary variables to an optimizer in the form of an integer linear program. The algorithm is outlined in pseudocode in Appendix A.  All orbits are discretized such that successive position states are 0.015 TU $\approx$ 93 minutes apart.

The rest of this section first introduces the CR3BP dynamics and then introduces how each input into the ILP formulation is defined for the cislunar SDA architecture design problem.

\subsection{Circular Restricted 3 Body Problem}
To leverage this formulation in cislunar space, a gravitational model of the region is needed. The CR3BP offers a reasonably accurate model for the dynamics of a small satellite as it moves under the influence of two much larger bodies (in our case, the Earth and Moon).  The CR3BP assumes a point of negligible mass under the gravitational influence of two much larger point masses that have circular orbits about their common barycenter. The equations are given by,

\begin{equation}
\begin{aligned}
    \ddot{x} - 2\dot{y} &= \frac{dU}{dx}\\
    \ddot{y} + 2\dot{x} &= \frac{dU}{dy}\\
    \ddot{z} &= \frac{dU}{dz}
\end{aligned}
\end{equation}

where $U$ is an effective potential of the rotating system, given by
\begin{equation}
    U = \frac{x^2 + y^2}{2} + \frac{1-\mu}{r_1} +\frac{\mu}{r_2}
\end{equation}

Here, $\mu = \frac{m_2}{m_1 + m_2}$ is the mass ratio of the primaries with masses of $m_1$ and $m_2$, and $r_1$ and $r_2$ are the distances of the object to each primary. 

\begin{equation}
    \begin{aligned}
        r_1 &= \sqrt{(x+\mu)^2 + y^2 + z^2}\\
        r_2 &= \sqrt{(x-(1-\mu))^2 + y^2 + z^2}
    \end{aligned}
\end{equation}

This formulation of the CR3BP uses normalized mass, distance, and time units. It is important to note that the CR3BP admits periodic orbits as solutions, which is convenient for the ILP formulation developed above. Furthermore, the CR3BP admits 5 stationary points, called \textit{Lagrange} points, commonly referred to as L1 through L5.

For certain accessibility measures such as the apparent magnitude of the object, the position of the sun is required in addition to that of the observer and target. In the rotating frame of the CR3BP, the sun's motion is clockwise with a period equal to the synodic period of the Earth-Moon system. A simplified model of the sun's motion is a coplanar clockwise circular orbit of radius 1 AU around the Earth-Moon barycenter. Note that the sun is only used as an illumination source in this model. We do not consider the sun's gravitational perturbation on the system. Parameters for the dynamics of this system are gathered in Appendix B.

\subsection{Target Point Set $\mathcal{J}$}\label{subsec4}
The set of target points can be any spatio-temporally varying points, discretized in space and time. In this paper, we analyze two cases as application examples. 

\subsubsection{Case 1: Surveillance of Objects in Translunar Orbit }\label{subsubsec2}
\begin{figure}[h!]
	\centering\includegraphics[width=0.8\linewidth]{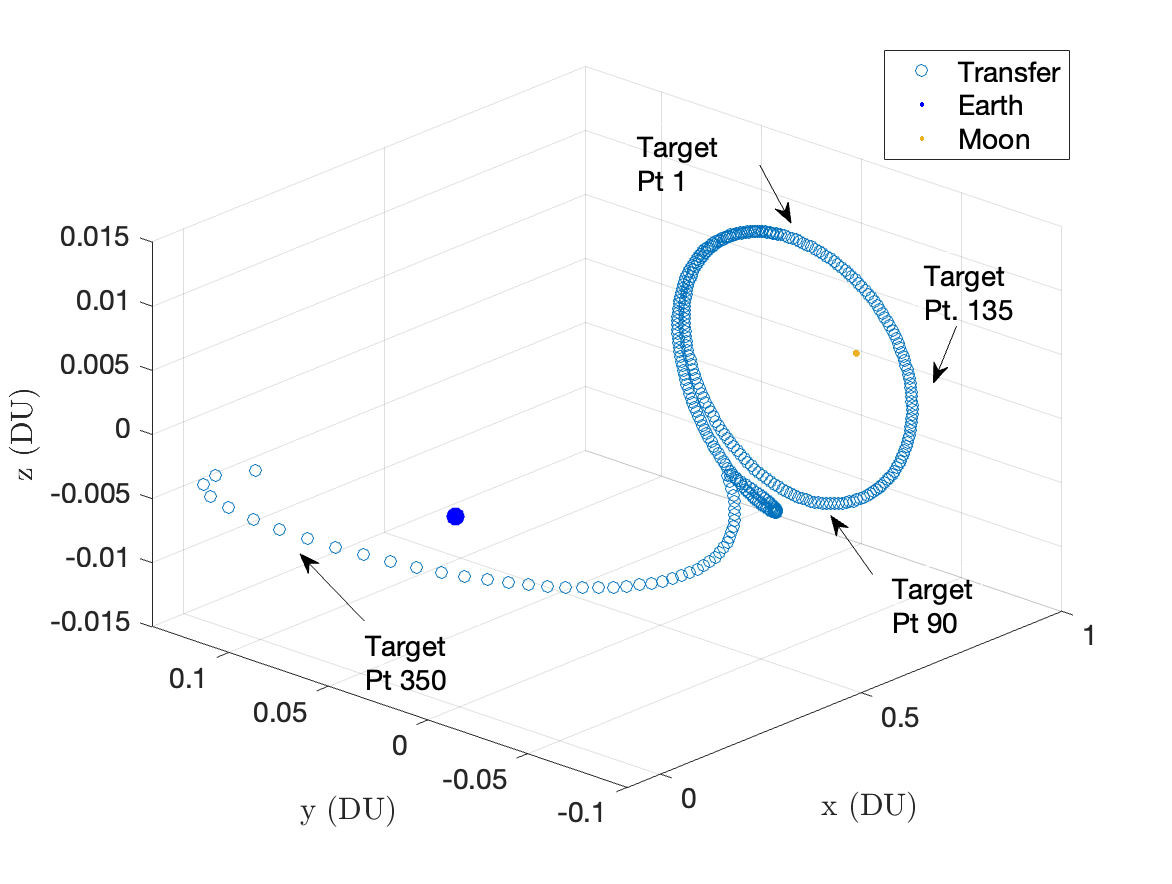}
	\caption{L1 to GEO transfer}
	\label{fig: 3}
\end{figure}

As traffic between the Earth and the Moon is anticipated to increase, an attractive choice is a transfer orbit from L1 to geosynchronous Earth orbit (GEO). Figure \ref{fig: 3} illustrates this transfer. It consists of a manifold insertion from an L1 Halo Orbit followed by a high impulse manuever into a connecting orbit to GEO. The transfer takes $\approx$ 23.12 days.
The discretized position history (i.e. the blue circles in the figure) are target points comprising the set $\mathcal{J}$. Note that although we define $\mathcal{J}$ as a \textit{set} of points, it is actually a \textit{sequence}. Target point 1 is followed by target point 2 and so on, tracking the position of a hypothetical target on this path from L1 to GEO. If we want to cover multiple target departure time possibilities (referred to as \textit{departure time windows}), we would have multiple such target trajectory candidate time sequences over this orbit, one corresponding to each departure time possibility. This transfer is referred to as ``Halo2GEO'' in this paper.

\subsubsection{Case 2: Simultaneous Surveillance of Multiple Orbits}\label{subsubsec3}
\begin{figure}[H]
	\centering\includegraphics[width=\linewidth]{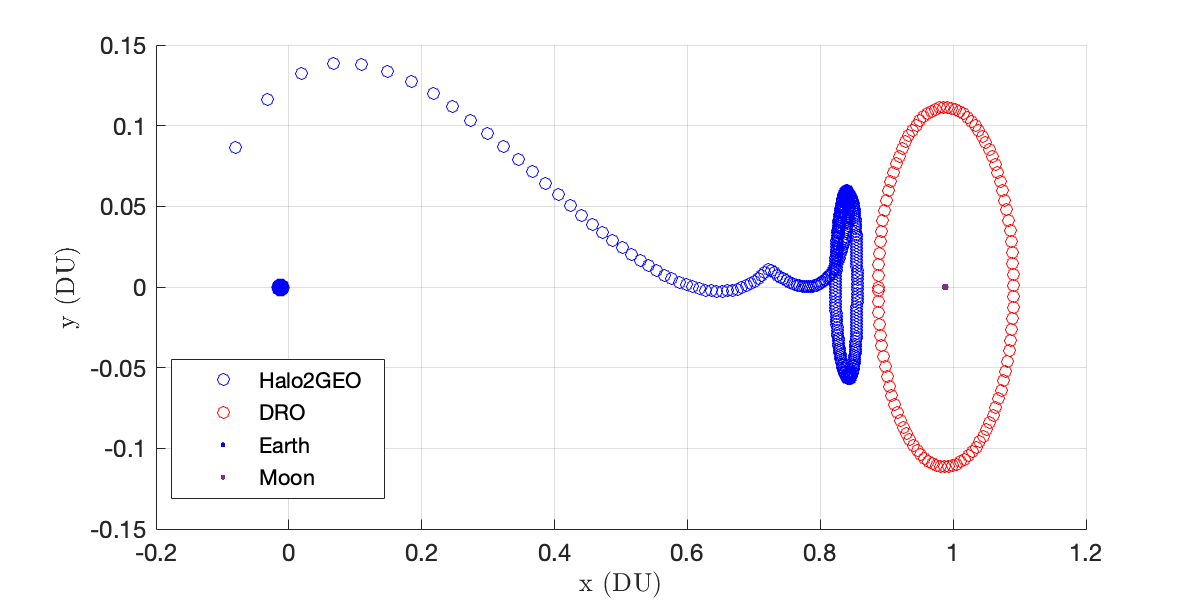}
	\caption{DRO and Halo2GEO Target Points}
	\label{fig:4}
\end{figure}

To showcase the power of the formulation, we consider a simultaneous coverage task, realized as two separate targets on separate trajectories. The first is the Halo to GEO transfer from Case 1. The second is a distant retrograde orbit (DRO) around the Moon, which can, for example, represent critical lunar communications infrastructure in orbit. The target points are shown in Figure \ref{fig:4}. This is the target point set $\mathcal{J}$ for Case 2. It should be noted that the Halo2GEO transfer and the DRO orbital movement happen simultaneously. As the target moves from the Halo orbit to GEO, another target is moving on the DRO at the same time, hence the \textit{simultaneous} coverage task.
We optimize first for Halo2GEO only (Case 1), then DRO only, and finally for both Halo2GEO and DRO considered concurrently to investigate how the optimizer removes any  redundancies  arising from optimizing constellations separately.

\subsection{Coverage Requirement $\boldsymbol{f}$}\label{subsec4}
Given the target points in the set $\mathcal{J}$, We need to define the coverage requirement $\boldsymbol{f}$ for each target point. Namely, we seek a set $\mathcal{F} = \{ \boldsymbol{f}_1, \boldsymbol{f}_2, \cdots, \boldsymbol{f}_{\lvert\mathcal{J}\rvert} \}$ of coverage requirements, one for each target point in the set $\mathcal{J}$. In this paper, we assume that we aim to track a target departing on an orbit (e.g., Halo2GEO orbit) but at different possible departure times. The target passes through each target point consecutively.  Thus, the coverage requirement for target point $j+1$ should be that of target point $j$, but shifted (in fact, circularly shifted due to the periodicity of observer orbits) forward by one index. Mathematically,

\begin{equation}
    f_{j+1}[i+1] = f_j[i \bmod L] \label{eqn:cov_req}
\end{equation} 

Where $[\ ]$ symbolizes indexing into the coverage requirement array $\boldsymbol{f}_j$. This recursive relation is convenient because we are only required to define the coverage requirement $\boldsymbol{f}_1$, with subsequent $\boldsymbol{f}_j$ defined by Eqn. \ref{eqn:cov_req}. An example set $\mathcal{F} = \{ \boldsymbol{f}_1, \boldsymbol{f}_2, \boldsymbol{f}_3 \}$ is illustrated below in Figure \ref{fig:6}.
\begin{figure}[H]
	\centering\includegraphics[width=\linewidth]{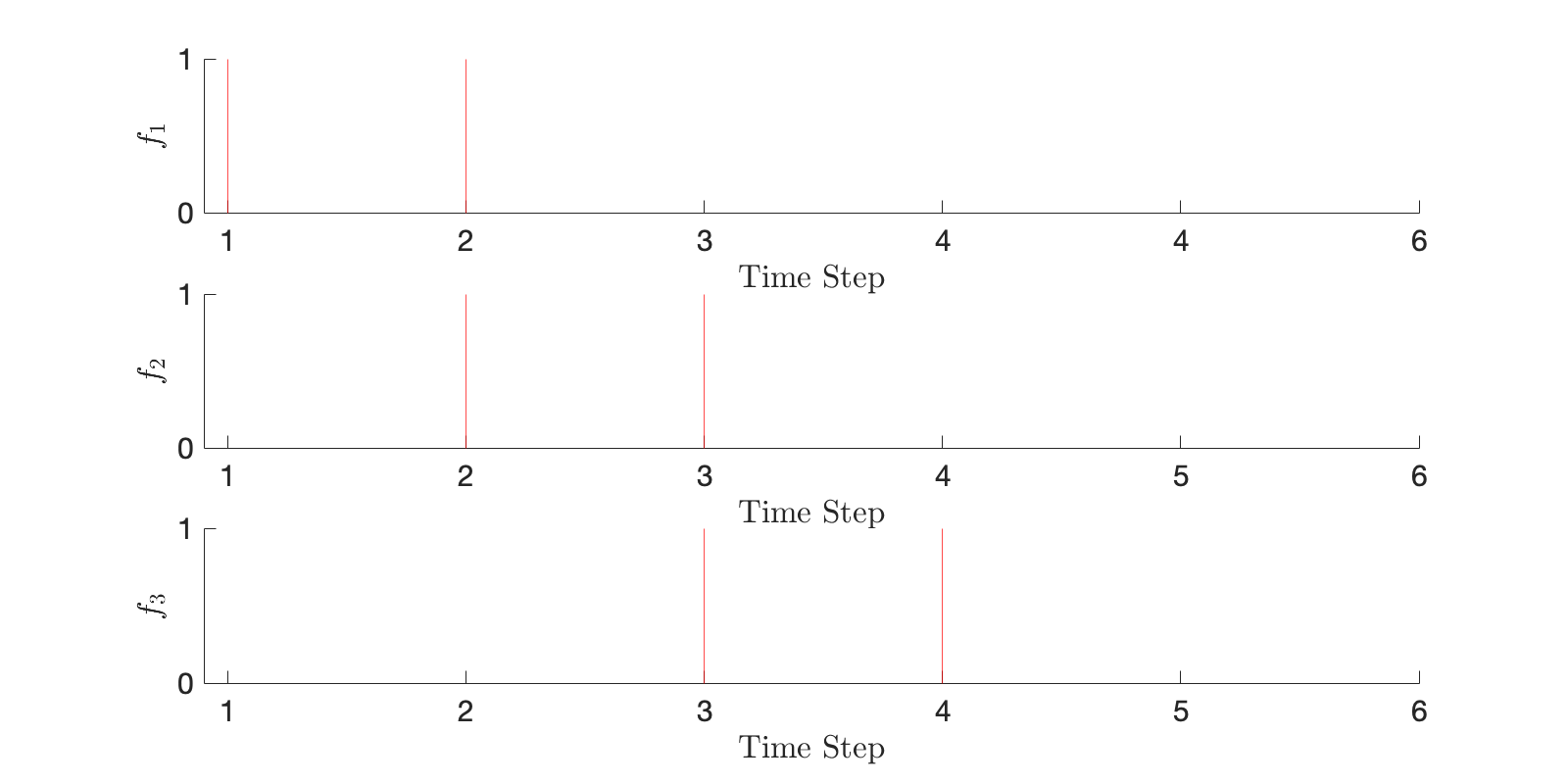}
	\caption{Illustration of Equation 1: an example  where each target point requires two sequential observations. The requirements shift right as a hypothetical target moves through each target point. }
	\label{fig:6}
\end{figure}

\subsection{Observer Orbit Set $\mathcal{Z}$}\label{subsec4}
There are a variety of periodic orbits in the CR3BP. Of these, the most attractive options offer close approaches to L1 and GEO while moving through the majority of cislunar space. We choose a set of 6 orbits inspired by Gupta et al.\cite{guptaconstellation}, which are plotted in Figure \ref{fig:5} below. 

\begin{figure}[H]
	\centering\includegraphics[width=\textwidth]{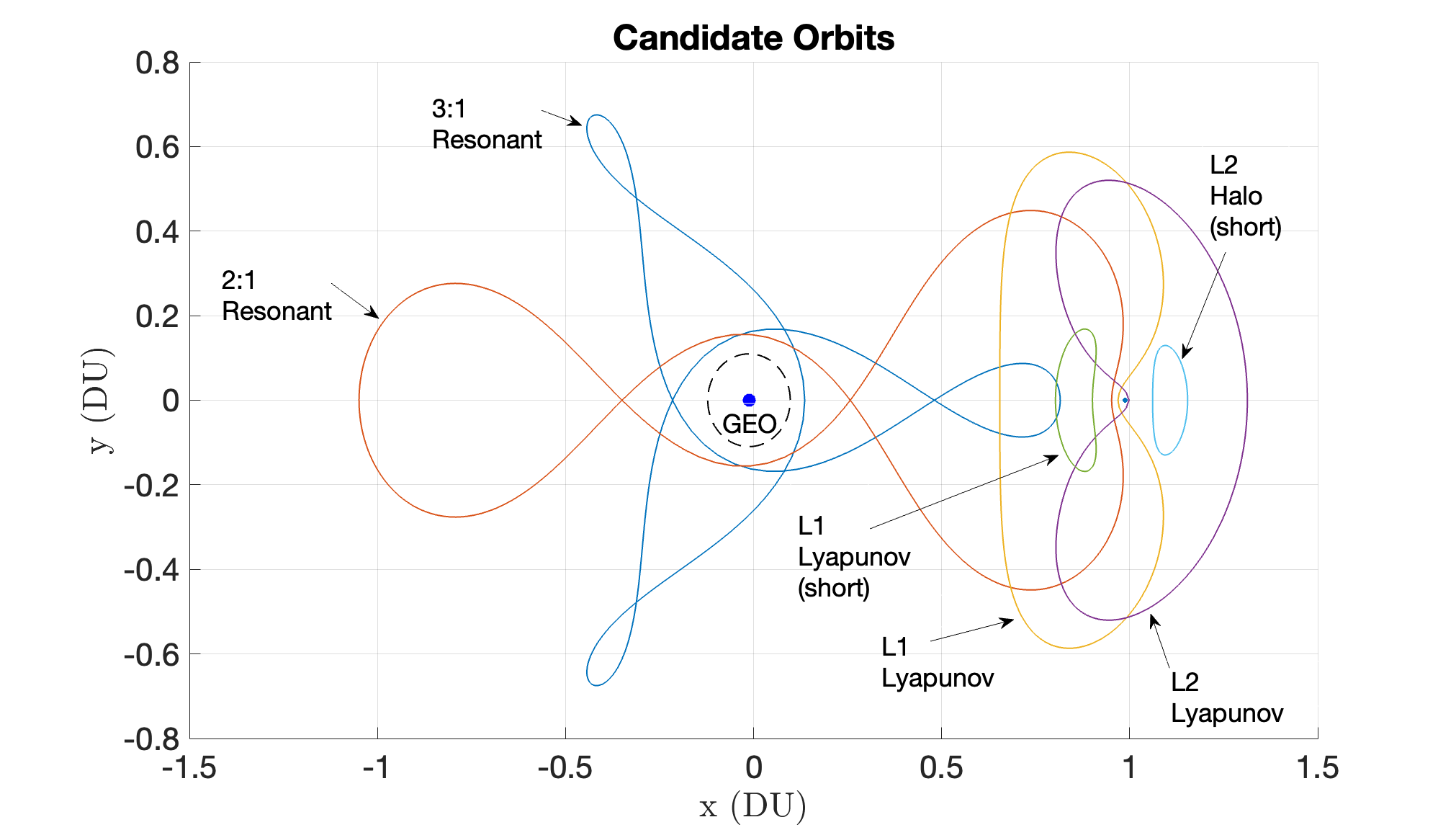}
	\caption{Candidate Orbits.  Orbits have the same period in the rotating frame to allow an optimized constellation's configuration to repeat. Thus, the optimized constellation's coverage guarantees can extend beyond the simulation time.}
	\label{fig:5}
\end{figure}
The initial conditions for these orbits are shown in Appendix B. All orbits were generated via a single-shooting differential correction algorithm with initial guesses sampled from data provided by NASA JPL \cite{48975_2018} and Broucke \cite{1968port.book.....B}. All orbits have period $T$ = 6.45 TU $\approx$ 28.008 days in the rotating frame; note that one of the considered L1 Lyapunov orbits and the considered L2 Halo orbit have a fundamental period of $\frac{T}{2}$ = 3.225 TU which means they are also $T$-periodic. They have been labeled with ``short" qualifiers. In general, there is limited diversity in orbits with a period exactly equal to the synodic period of the Earth-Moon system (approximately 29.5 days).  As a consequence, while the orbits themselves repeat, the illumination conditions will not necessarily. As such, the results of the formalism apply on a case-by-case basis depending on the sun's initial phase (in radians) relative to the Earth-Moon axis, denoted as $\phi_0$. We stress that the candidate orbits are chosen to have periods as close as possible to the synodic period and are still resonant in the CR3BP sense. This ensures the difference in illumination conditions is small after $t = 6.45$ TU compared to a randomly chosen constellation period. We aimed to minimize this difference while preserving candidate orbit diversity. In this work, we choose $\phi_0$ = 0 for optimization (Earth-Moon-Sun are collinear at simulation start). We investigate the robustness of the constellation to different initial phase angles later in the results in section \ref{results:rob_subphase}. 

\subsection{Accessibility Metric and Threshold}\label{subsec5}
To quantify the observability of a target from an observer's position, we utilize the target point's apparent magnitude \cite{1974STIN...7512024K, vendl2021cislunar}. Consider an observer at position $\boldsymbol{r}$, a spherical target with diameter $d$ at position $\boldsymbol{p}$, and the sun at position $\boldsymbol{s}$. Define the \textit{solar phase angle} $\psi$,

\begin{equation}
    \cos{\psi} = \frac{(\boldsymbol{p} - \boldsymbol{r}) \cdot (\boldsymbol{p} - \boldsymbol{s})}{\|\boldsymbol{p} - \boldsymbol{r} \| \| \boldsymbol{p} - \boldsymbol{s}\|}
\end{equation}

It should be noted that $\psi$ is the familiar angle used in the dot product formula between two vectors, taking on values between 0 and $\pi$.
The diffuse phase angle function  \cite{1974STIN...7512024K}  can be written as a function of $\psi$,

\begin{equation}
    p_{\text{diff}}(\psi) = \frac{2}{3\pi}(\sin{\psi} + (\pi - \psi) \cos{\psi})
\end{equation}

On the interval (0, $\pi$), $p_{\text{diff}}$ is a strictly decreasing function. The apparent magnitude of a target \cite{1974STIN...7512024K}  is given by,

\begin{equation}   
    g(\boldsymbol{p}, \boldsymbol{r}, \boldsymbol{s}) \triangleq  m_{\text{sun}} - 2.5 * \log_{10}  \Big[ \frac{d^2}{\zeta^2} *(\frac{a_{\text{spec}}}{4} + a_{\text{diff}}*p_{\text{diff}}(\psi)) \Big]
\end{equation}
where $\zeta = \|\boldsymbol{p} - \boldsymbol{r} \|$ is the distance between the observer and target, and $a_{\text{spec}}, a_{\text{diff}} $ are the specular and diffuse reflection coefficients of said target, respectively. Values can be found in Appendix B. It should be noted that the solar phase angle $\psi$ is distinct from the initial phase angle $\phi_0$ described in section \ref{subsec4} above. $\phi_0$ describes the initial orientation of the Sun relative to the Earth-Moon axis. $\psi$ describes the orientation of an observer, a target point, and the Sun. The value of $g$ increases with $\psi$, meaning solar phase angles closer to $\pi$ result in an unfavorably dim target. Every unit increase in apparent magnitude corresponds to a 2.5 times dimmer target. In addition to this, we consider the possibility that the Earth and Moon may occlude the target. If this is the case, the accessibility measure is set to $\infty$, indicating an invisible target point. In practice, the measure will be set to the largest floating point number (often denoted REALMAX) available on the machine. We choose a threshold value of $M=17$ as our access measure cutoff. $g$ and $M$ define our accessibility metric and threshold, respectively.

\section{Results}\label{sec4}
This section discusses the results from each case discussed in Section \ref{sec3}. Case 1 analyzes the optimal constellation for monitoring the Halo2GEO orbit with multiple departure windows, while Case 2 analyzes the value of using a common constellation to monitor multiple target orbits: Halo2GEO + DRO.
\subsection{Case 1: Surveillance of Objects in Translunar Orbit}\label{subsubsec4}
We first focus on Case 1, with the Halo2GEO orbit. We will investigate the optimal satellite constellation to cover multiple departure time windows (i.e., multiple trajectory candidate time sequences), and investigate how this number of departure windows affects the satellite constellation size. More departure windows reflect an increasing uncertainty about when the target will enter its transfer orbit.  With more departure windows, the constellation is optimized for robustness against target departure delays. However, this will come at the cost of increasing constellation size. Thus, through this analysis, we can show the value of knowing precisely the departure time (or the cost of design for uncertainties in the departure time). 

\begin{figure}[H]
    
  \centering
  \begin{minipage}[b]{0.45\textwidth}
    \includegraphics[width=\textwidth]{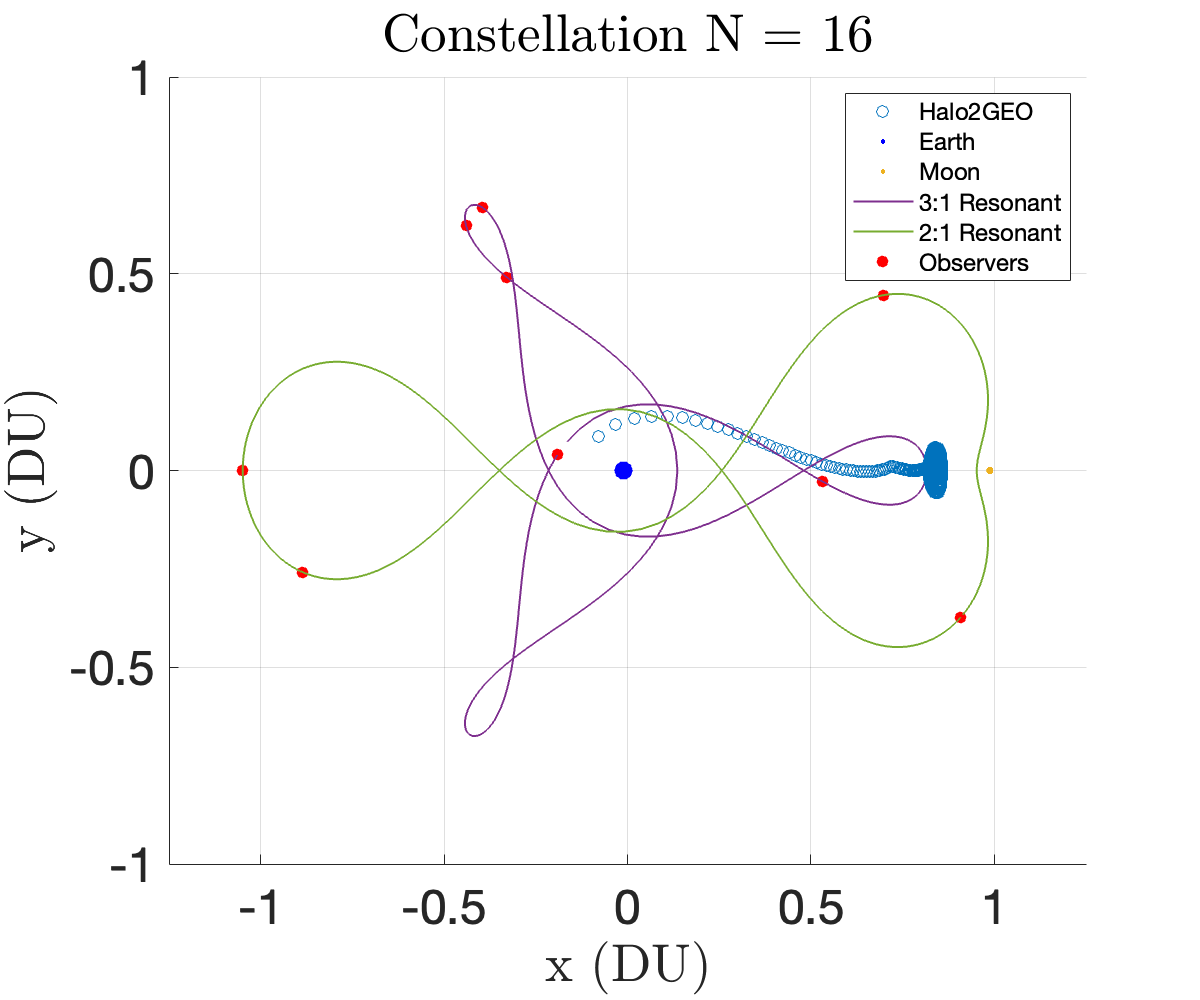}
    \caption{Optimal Constellation for 16 Departure Windows.  Observers are placed on 2:1 and 3:1 resonant orbits, likely due to close approaches to the Earth and Moon.}
    \label{fig:7}
  \end{minipage}
  \hfill
  \begin{minipage}[b]{0.45\textwidth}
    \includegraphics[width=\textwidth]{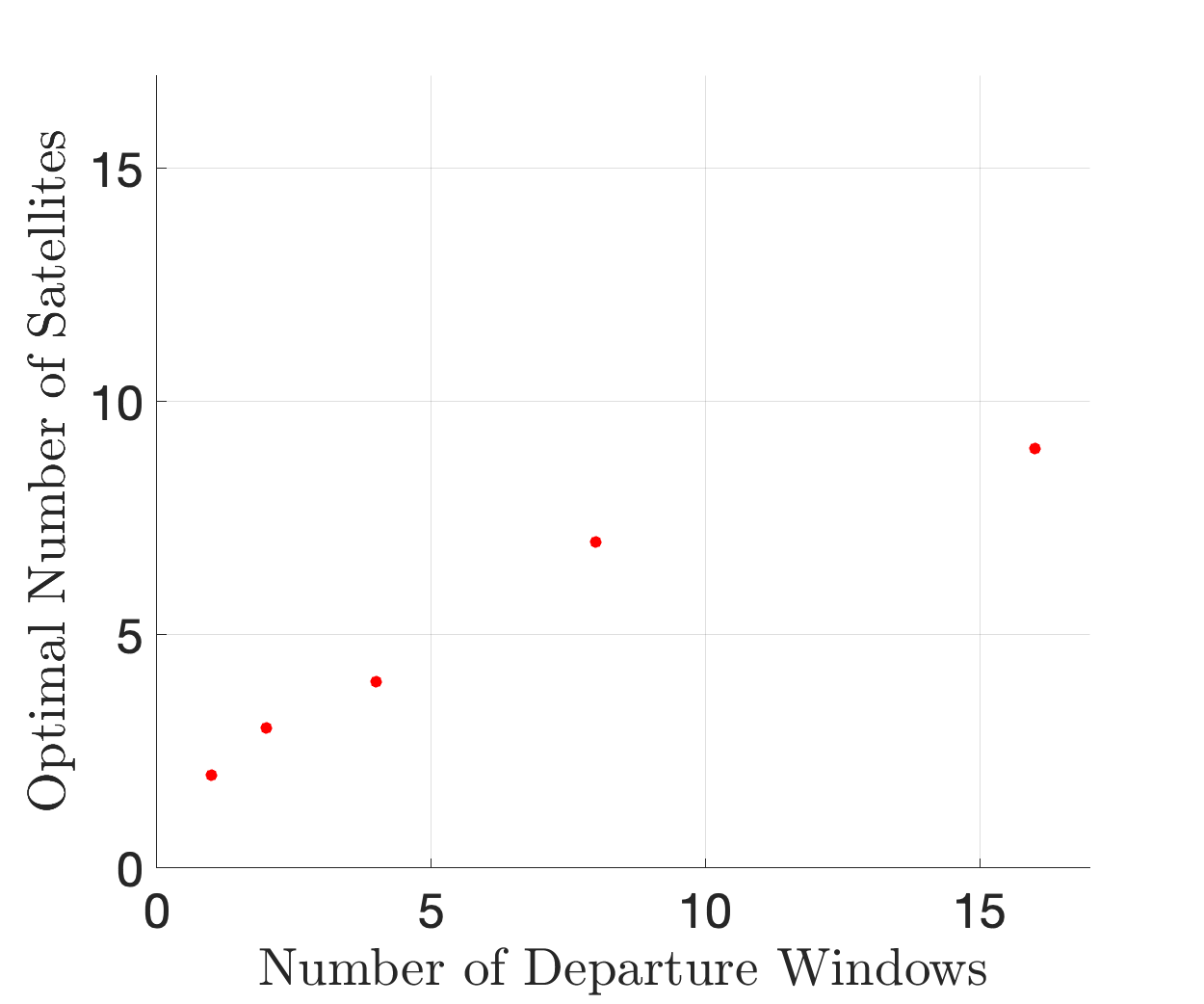}
    \caption{Constellation Size vs. Number of Departure Windows.  The optimal number of satellites grows sublinearly with the number of departure windows, likely because at some finite constellation size, coverage ability saturates.}
    \label{fig:8}
  \end{minipage}
\end{figure}

Denote the number of departure windows by $N$.  We specifically investigate $N \in\{1, 2, 4, 8, 16\}$, where the departure windows are uniformly separated amongst the $L$ time steps The powers of two ensure that the indices where $f_j = 1$ for $k$ departure windows are a subset of the indices where $f_j = 1$ for $2k$ departure windows, all while keeping the windows uniformly distributed in the interval $[0,  L]$. An algorithm for generating coverage requirement $\boldsymbol{f}_1$ is gathered in Appendix C. With $\boldsymbol{f}_1$, one can generate the entire coverage requirement using Equation \ref{eqn:cov_req}.  Figure \ref{fig:7} illustrates the result of the optimizer for $N = 16$ departure windows. The constellation consists of 9 satellites distributed on 3:1 and 2:1 resonant orbits, likely because of their close approaches to the moon, where a significant portion of the spatio-temporal coverage requirement is located. Figure \ref{fig:18} in Appendix D illustrates the results of the optimizer for every element $N \in\{1, 2, 4, 8, 16\}$. Figure \ref{fig:8} illustrates the optimal number of satellites versus the requested number of departure windows. The constellation size grows with the coverage demand. This is intuitive since coverage potential grows with the number of observers in the constellation. 

As the coverage requirement becomes continuous, the optimal number of satellites is expected to plateau. At some finite constellation size, it is expected that coverage ability saturates. Significantly stronger computational power is necessary to analyze more departure windows. A continuous coverage requirement, i.e. $\boldsymbol{f} = \boldsymbol{1}$, would require parallel computation for reasonable evaluation time.  We hypothesize that more complex coverage requirements lead to smaller solution sets among the large design space. Furthermore, members of the solution set will satisfy coverage requirements more tightly.

\begin{figure}[H]
    
  \centering
  \begin{minipage}[b]{0.45\textwidth}
    \includegraphics[width=\textwidth]{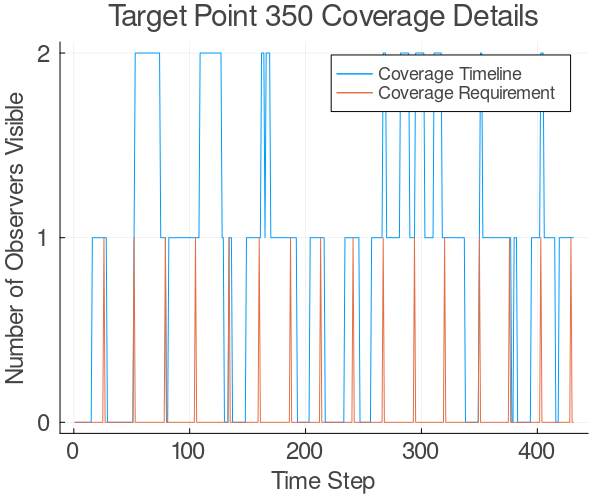}
    \caption{Coverage Details for Target Point 350 in a constellation optimized for 16 Departure Windows. Red impulses show the coverage requirement. The blue curve shows the coverage timeline.}
    \label{fig:a}
  \end{minipage}
  \hfill
  \begin{minipage}[b]{0.45\textwidth}
    \includegraphics[width=\textwidth]{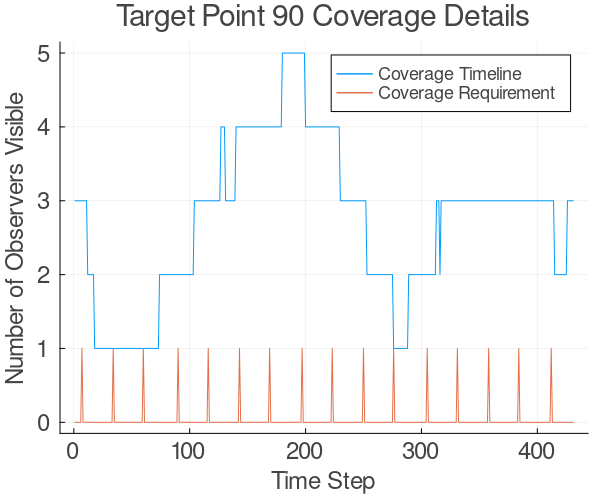}
    \caption{Coverage Details for Target Point 90 in a constellation optimized for 16 Departure Windows. Coverage is relatively more saturated.}
    \label{fig:b}
  \end{minipage}
  
\end{figure}

It is also interesting to investigate which \textit{portions} of the coverage requirement the solution satisfies more tightly. Informally, this reveals which target points are more difficult to observe. To this aim, we illustrate the saturation of coverage at different target points for the constellation shown in Figure \ref{fig:7}. Specifically, we show the coverage timelines and coverage requirements of two target points (target point 350 and 90) in Figure \ref{fig:a} and Figure \ref{fig:b}, respectively. We see that coverage is sparse (i.e the constraint $\boldsymbol{Vx} \geq \boldsymbol{f}$ is satisfied more tightly) for target point 350 compared to target point 90. Target points closer to GEO are more difficult to monitor compared to those near the Earth-Moon L1. A more complete visualization is presented later in section \ref{subsec5}. 

\begin{figure}[H]
    
  \centering
  \begin{minipage}[b]{0.45\textwidth}
    \includegraphics[width=\textwidth]{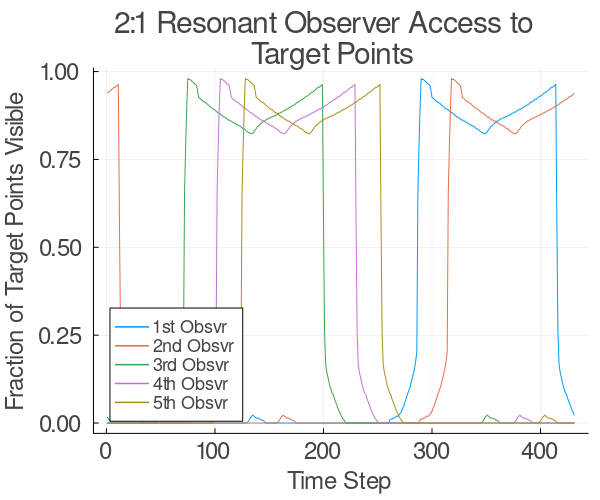}
    \caption{Time evolution of each 2:1 resonant observer's access to all target points. }
    \label{fig:c}
  \end{minipage}
  \hfill
  \begin{minipage}[b]{0.45\textwidth}
    \includegraphics[width=\textwidth]{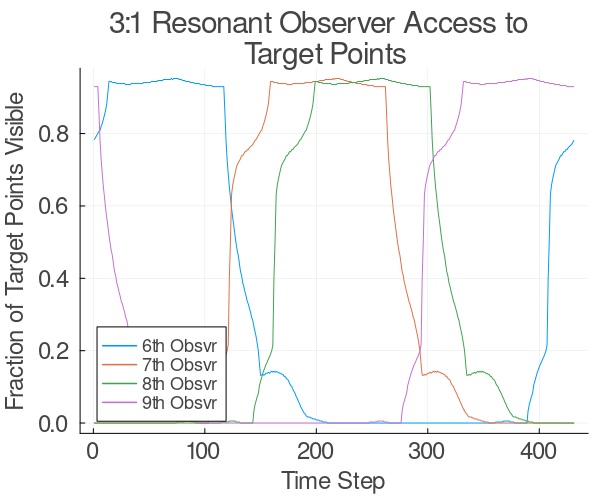}
    \caption{Time evolution of each 3:1 resonant observer's access to all target points. }
    \label{fig:d}
  \end{minipage}
  
\end{figure}

Above we presented how, for a chosen target point, the number of available observers fluctuates across the simulation time. It is also interesting to investigate the converse: for a chosen observer, how the number of target points in view fluctuates across the simulation time. For the constellation shown in Figure \ref{fig:7}, observers are placed on 2:1 and 3:1 resonant orbits. Figure \ref{fig:c} shows the fraction of available target points for each observer on the 2:1 resonant orbit as a function of time. Figure \ref{fig:d} shows the same quantity but for each observer on the 3:1 resonant orbit. Intuitively, in each plot, the curve for each observer is identical but shifted across time to reflect the relative phasing between observers. The figures show that each observer cannot see all the target points at a given time, reflecting the need for the satellites in the constellation to work \textit{together} to satisfy coverage demand.

\begin{figure}[H]
    
  \centering
  \begin{minipage}[b]{0.45\textwidth}
    \includegraphics[width=\textwidth]{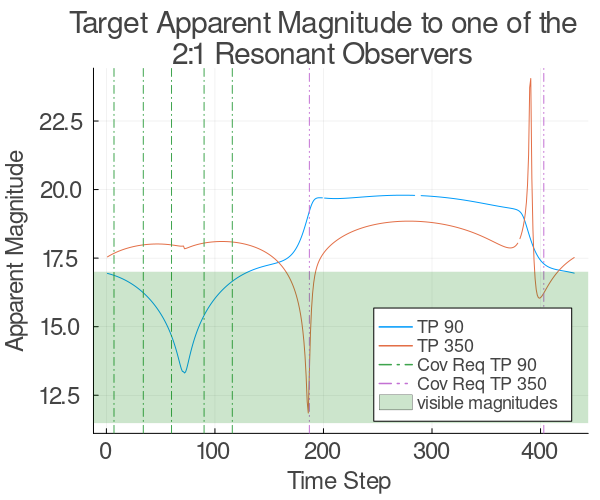}
    \caption{Apparent magnitude of two target points to a 2:1 resonant observer. Discontinuities are due to target occlusion by the Earth and/or Moon. \textbf{Note}: only parts of the coverage requirement the selected observer helps satisfy are plotted (vertical dashed lines).}
    \label{fig:e}
  \end{minipage}
  \hfill
  \begin{minipage}[b]{0.45\textwidth}
    \includegraphics[width=\textwidth]{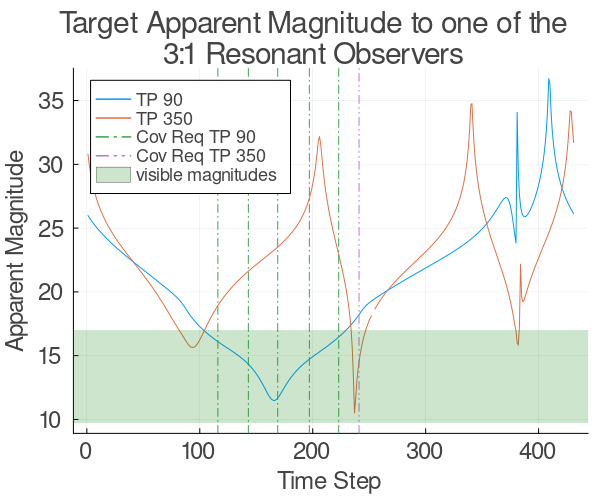}
    \caption{Apparent magnitude of two target points to a 3:1 resonant observer. Discontinuities are due to target occlusion by the Earth and/or Moon. \textbf{Note}: only parts of the coverage requirement the selected observer helps satisfy are plotted (vertical dashed lines).}
    \label{fig:f}
  \end{minipage}
  
\end{figure}

The reason why multi-observer collaboration is required for the constellation in Figure \ref{fig:7} can be explained by the restrictive solar illumination conditions for target observation. Figures \ref{fig:e} and \ref{fig:f} show the time evolution of the apparent magnitude of two target points for one of the 2:1 and 3:1 resonant observers, respectively. Target points are visible when the apparent magnitude is below the threshold (in the shaded green region). Note that both target points are unobservable for a significant portion of the simulation period (the majority of blue and red curves are not in the green region). Nevertheless, the selected observer plays the necessary role to satisfy the visibility requirements at the constellation level. To show that, the observations by the selected observer used to satisfy the coverage requirements of each target point are represented by the vertical dashed lines in the figures; we see that the apparent magnitude at these time steps for the corresponding target points is below the threshold, and thus the corresponding target points are visible. Note that only a subset of the requirement is satisfied by the selected observers, with the remaining requirement satisfied by the other observers.  

\subsection{Case 2: Simultaneous Surveillance of Multiple Orbits}\label{subsubsec5}

\begin{figure}[H]
  \centering
  \begin{minipage}[b]{0.45\textwidth}
    \includegraphics[width=\textwidth]{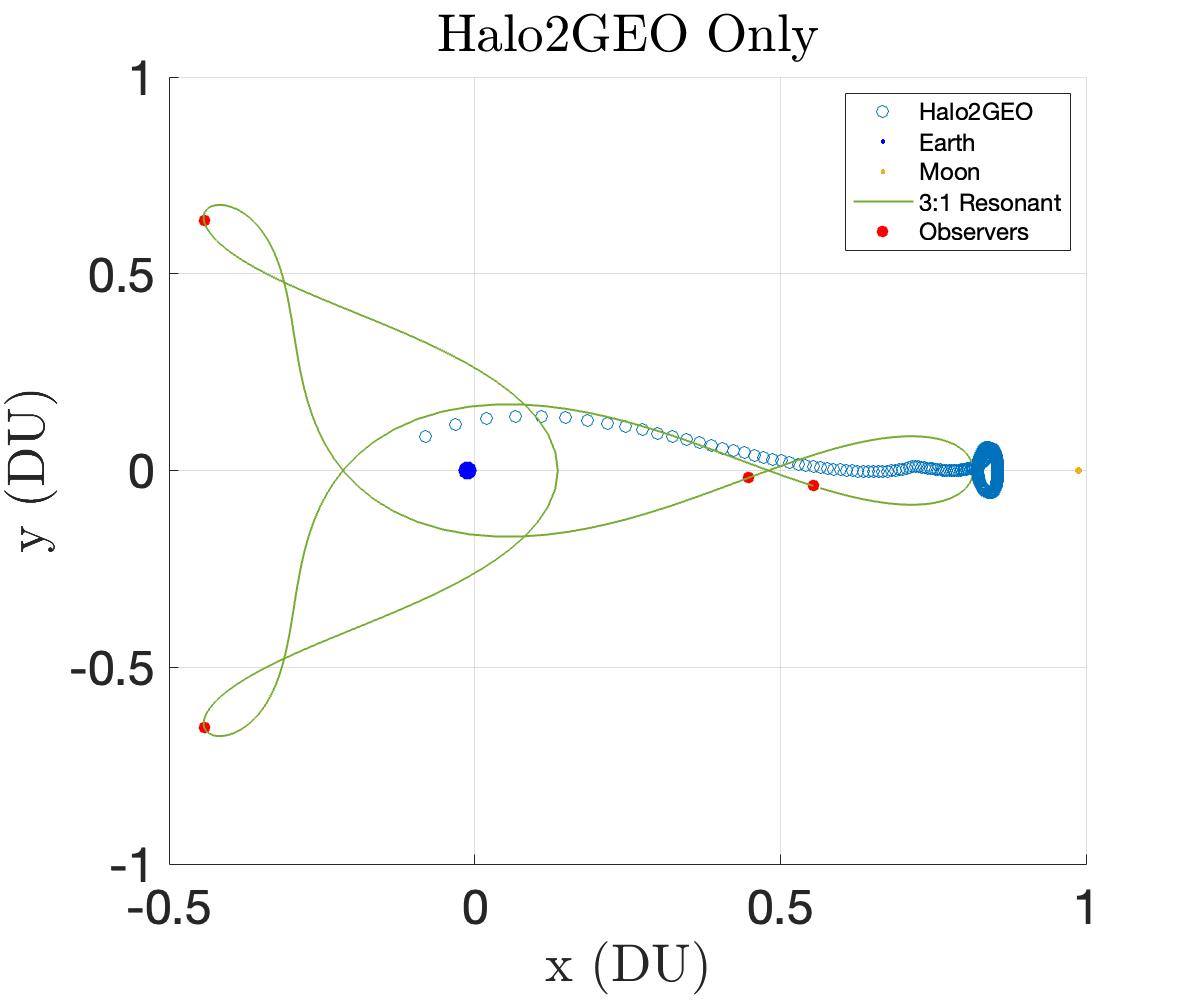}
    \caption{Optimal Constellation for Halo2GEO Transfer (4 satellites)}
    \label{fig:9}
  \end{minipage}
  \hfill
  \begin{minipage}[b]{0.45\textwidth}
    \includegraphics[width=\textwidth]{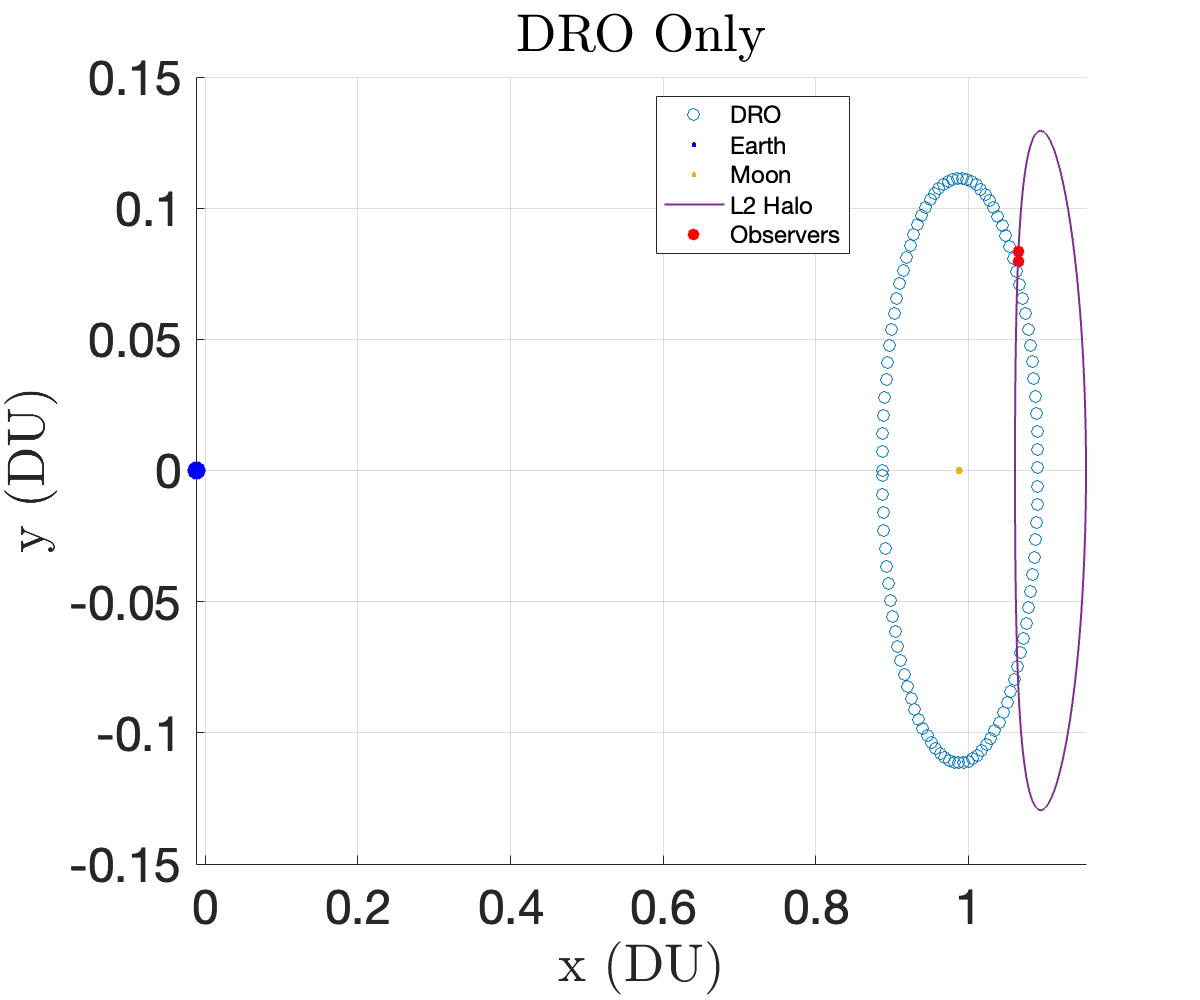}
    \caption{Optimal Constellation for Distant Retrograde Orbit (DRO) (2 satellites)}
    \label{fig:10}
  \end{minipage}
\end{figure}

\begin{figure}[H]
	\centering\includegraphics[width=0.8\textwidth]{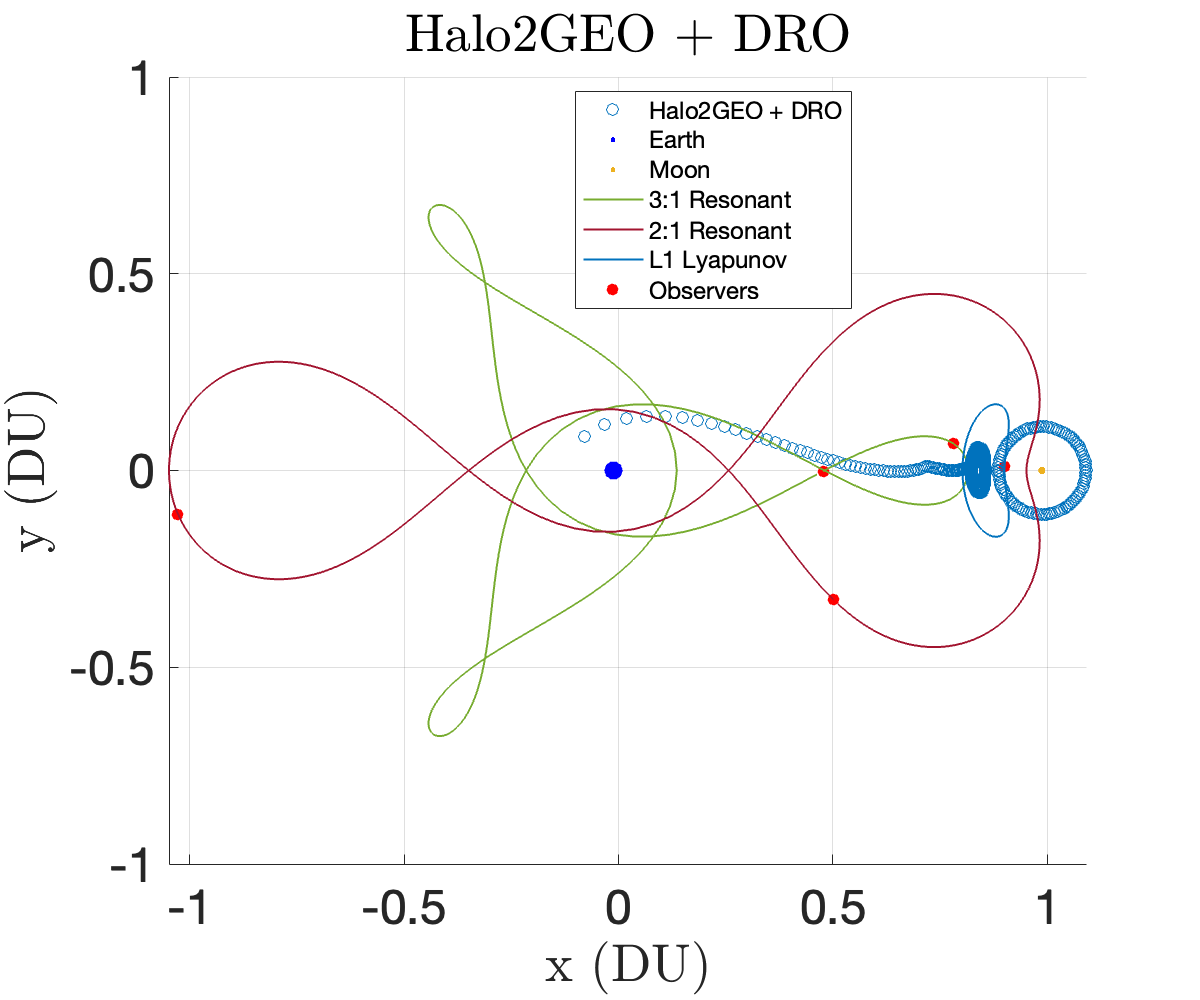}
	\caption{Concurrently Optimized Constellation (5 satellites).  When concurrently optimized, the constellation requires only 5 satellites instead of 6.}
	\label{fig:11}
\end{figure}

Suppose that the Halo2GEO transfer used N = 4 departure windows for its coverage requirement, while the DRO has only 1 departure window possibility but requires 2 satellites to always maintain custody. The second requirement can be stated mathematically as 
\begin{equation}
    \boldsymbol{f}_j[i] = 2  \delta_{ij}
\end{equation}
where $\delta_{ij} = 1$ if $i=j$ and 0 otherwise (i.e. the Kronecker delta).
Shown below are the separately optimized constellations for Halo2GEO only (Figure \ref{fig:9}) and DRO only (Figure \ref{fig:10}).
Monitoring the DRO alone requires 2 satellites situated close to each other to maintain the double custody requested by the requirement. Monitoring the Halo2GEO transfer requires 4 satellites and is equivalent to the ``Constellation N = 4'' constellation in Case 1. If both transfers are considered concurrently by the optimizer, it chooses the same set of orbits but phases the satellites differently, resulting in a constellation of only 5 satellites, shown in Figure \ref{fig:11}.

 In contrast to the small relative phasing between the two satellites in Figure \ref{fig:10}, we see the satellites are now distributed on separate orbits. The satellite on the L1 Lyapunov (short) orbit is efficient in satisfying parts of the coverage requirement for both the DRO and Halo2GEO transfer in addition to the assistance provided by satellites on the 2:1 and 3:1 resonant orbits. Reducing the constellation by 1 satellite provides an immediate launch cost relief in addition to a long-term cumulative maintenance cost relief. The proposed ILP formulation enables us to consider multiple targets concurrently, enabling the mission designer to identify and eliminate  redundancies  in constellations optimized separately for separate objectives.

\subsection{Robustness to Delays in Departure}\label{subsec5}
\begin{figure}[H]
  \centering
  \begin{minipage}[b]{0.45\textwidth}
    \includegraphics[width=\textwidth]{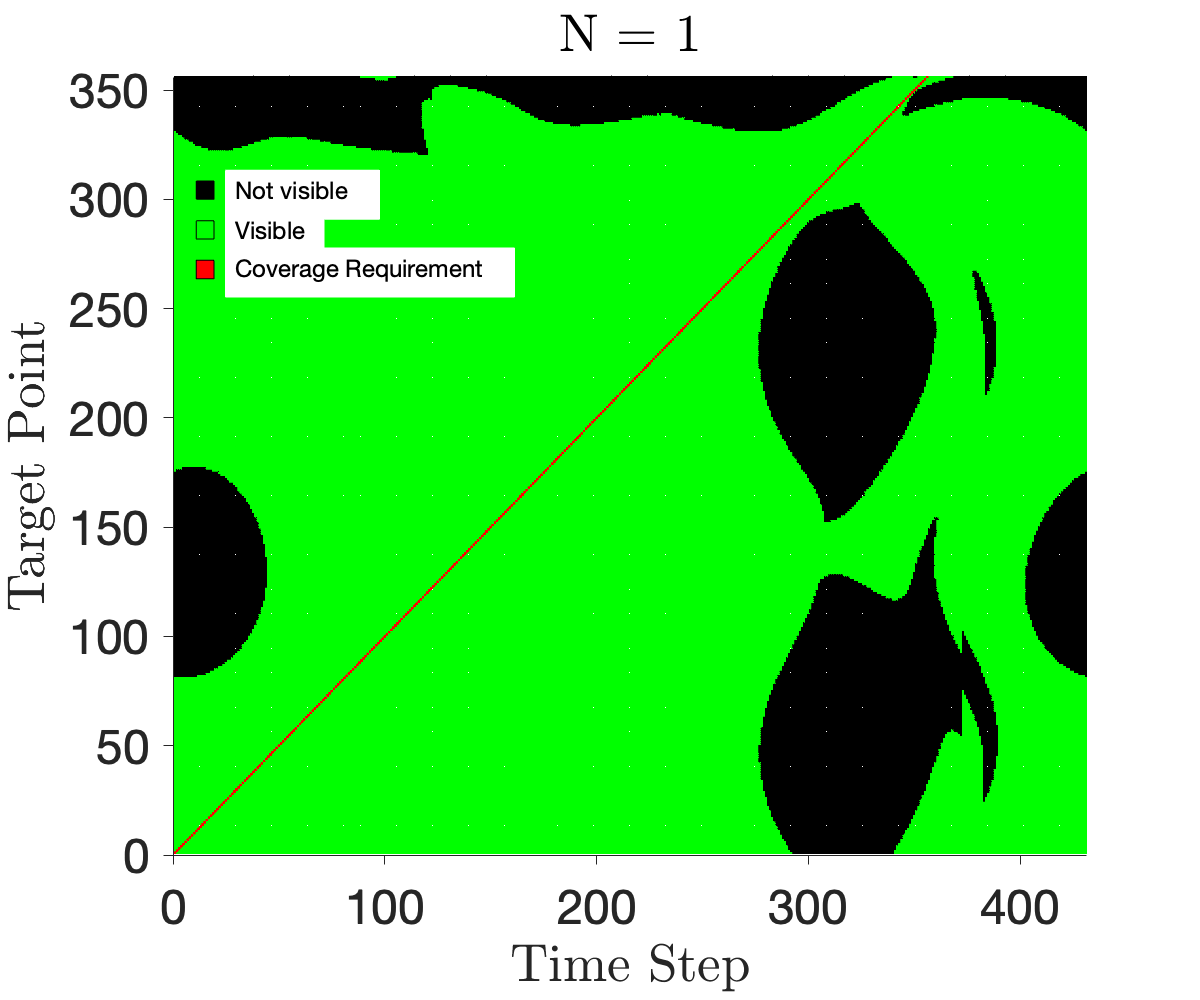}
    \caption{Coverage Robustness for 1 Departure Window.  Each point on the grid represents a \{time step, target point\} pair. If the point is green, there is at least one observer that can view the target at that time. If the point is black, the target is not visible to any observer at this time. Red points show the coverage requirement. }
    \label{fig:12}
  \end{minipage}
  \hfill
  \begin{minipage}[b]{0.45\textwidth}
    \includegraphics[width=\textwidth]{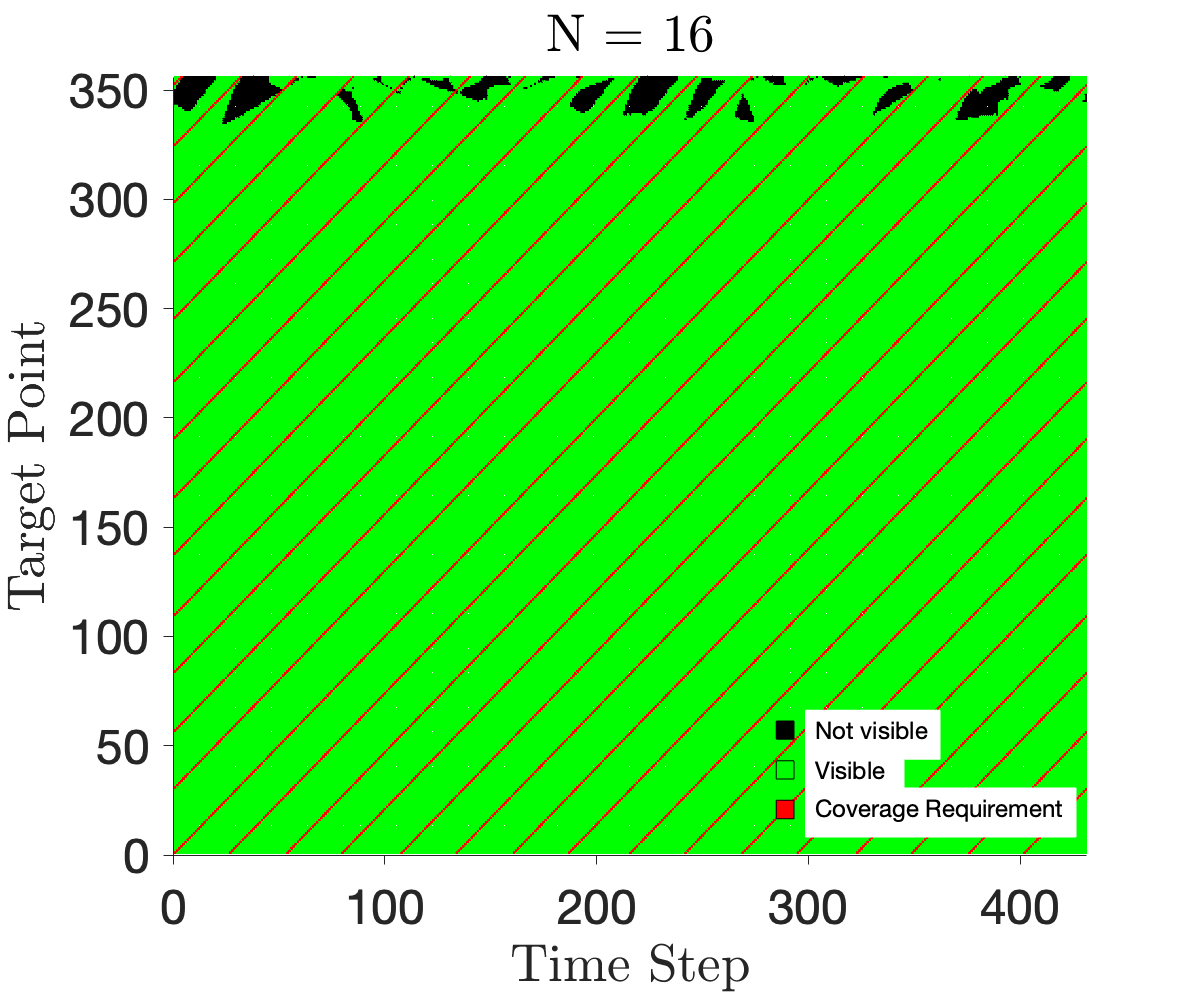}
    \caption{Coverage Robustness for 16 Departure Windows.  In addition to satisfying the coverage requirement, the constellation is able to observe target points at other times, as evidenced by the ratio of green to black in the grid.}
    \label{fig:13}
  \end{minipage}
\end{figure}

The constellation shows excellent robustness to delays in departure. The following results assume a Halo2GEO transfer only (i.e Case 1), the DRO is not part of $\mathcal{J}$. Figures \ref{fig:12} and \ref{fig:13} show binary indicator maps for constellations optimized for $N=1$ and $N=16$ departure windows. Figure \ref{fig:19} in Appendix D shows the maps for each $N \in \{1, 2, 4, 8, 16 \}$ departure windows. Areas shaded in green indicate a \{time step, target point\} pair which is observable by at least one satellite in the constellation. Areas shaded in black are not observable. The red lines (shown as diagonal streaks across Figures \ref{fig:12} and \ref{fig:13}) illustrate the spatio-temporal coverage requirement sent to the optimizer. These regions are always covered, as required by the formulation. We see that the majority of the optimizer's work attempts to cover the target points closer to GEO (target point No. 300 and onwards). As seen in Figure \ref{fig: 3}, these points are spaced further apart (the target is moving quickly in this region), naturally making it more difficult to observe. With increasing $N$, we make progress in closing the gap, with most of the binary heat map displaying green for $N = 16$.

\subsection{Robustness to Initial Sun Phase Angle, $\phi_0$} \label{results:rob_subphase}
\begin{figure}[H]
  \centering
    \begin{minipage}[b]{0.45\textwidth}
    \includegraphics[width=\textwidth]{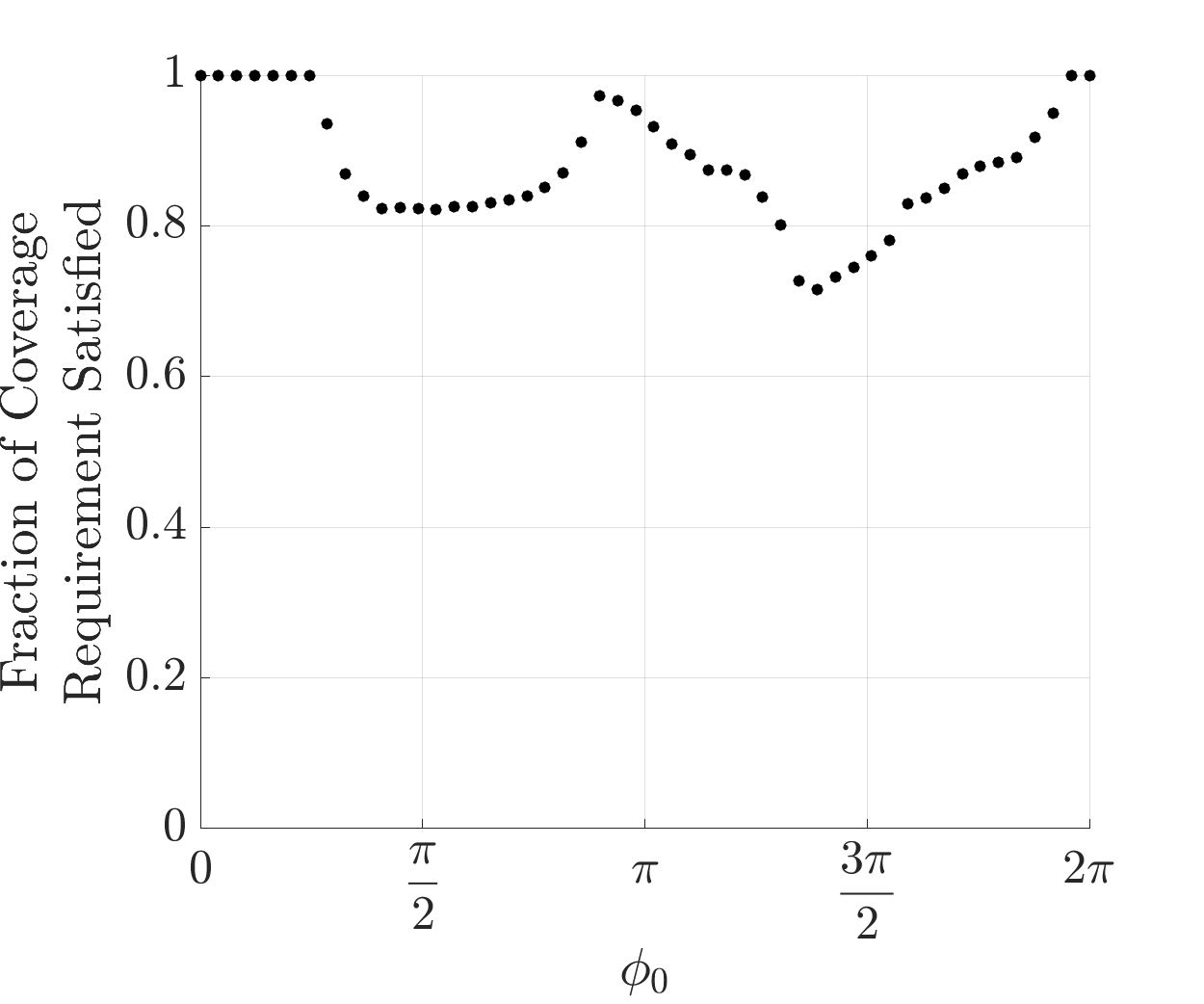}
    \caption{Coverage Requirement Satisfaction vs. $\phi_0$ for 2 Departure Windows.  Constellation exhibits impressive robustness. Although only optimized for one period, the worst-case performance of consequent periods is greater than 70\%.}
    \label{fig:14}
  \end{minipage}
    \hfill
  \begin{minipage}[b]{0.45\textwidth}
    \includegraphics[width=\textwidth]{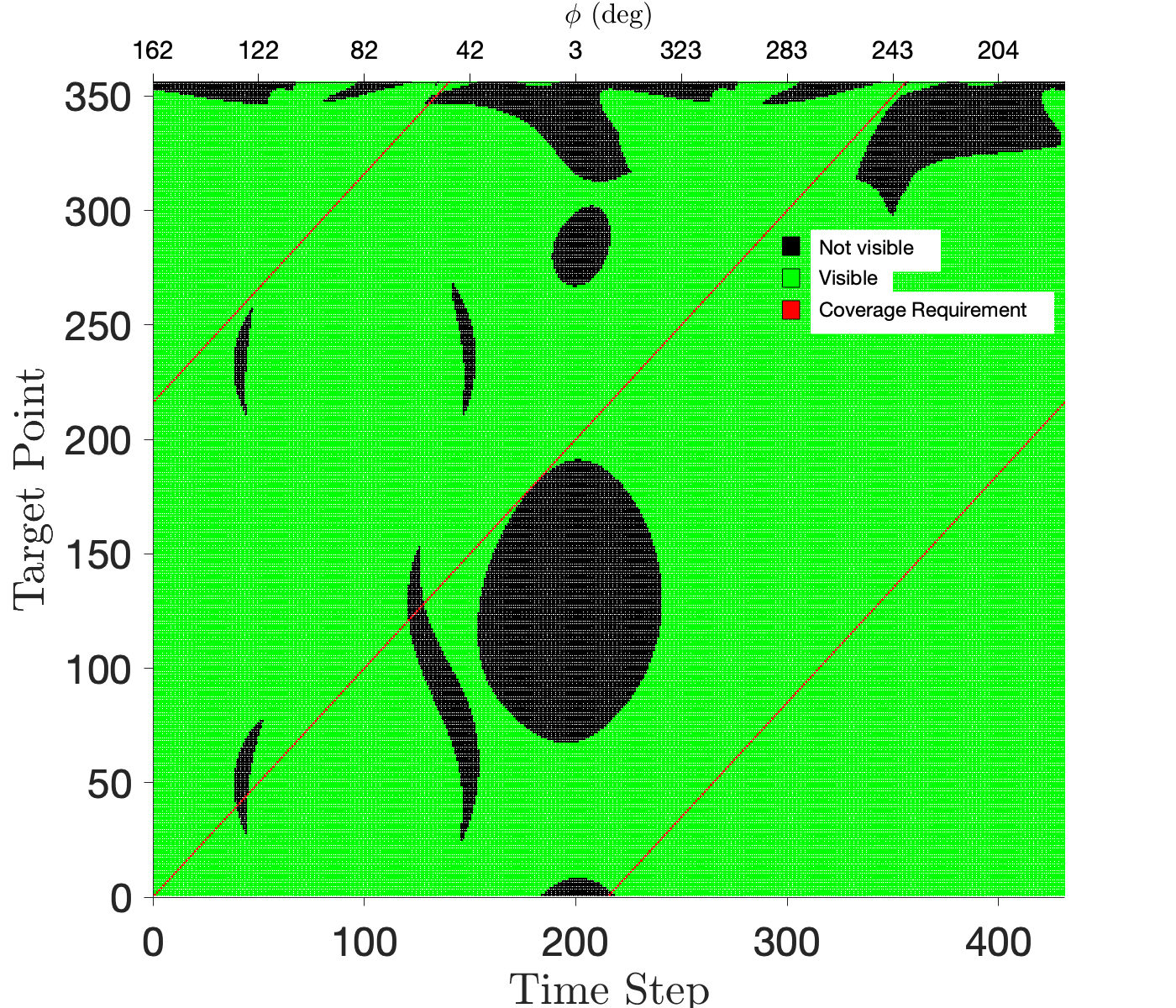}
    \caption{Coverage Robustness for 2 Departure Windows with $\phi_0 = 162^{\circ}$.  Constellation satisfies 97\% of coverage requirement.}
    \label{fig:15}
  \end{minipage}
\end{figure}

We acknowledged earlier that because our candidate orbits are not resonant with the Earth-Moon synodic period, after one simulation period, the illumination conditions will differ slightly at the start of the second epoch ($t = 6.45$ TU). The orbits repeat but now with a slightly different initial sun phase angle, $\phi_0$. Again, the following results assume a Halo2GEO transfer only (i.e., Case 1), the DRO is not part of $\mathcal{J}$. As discussed above, we chose an initial sun phase angle of 0 rad, meaning that at the start of the simulation, the Sun is on the positive x-axis at a distance of 1 AU. The constellations are optimized for this phase angle, but they show good robustness against various initial sun phase angles. We investigated numerous $\phi_0$ between 0 and 2$\pi$, given a constellation already optimized for $N=2$ departure windows (see Figure \ref{fig:18} in Appendix D for an illustration of the constellation). Figure \ref{fig:14} shows that despite changes in $\phi_0$, the constellation still satisfies a majority of the coverage requirement. As such, the constellation's effectiveness extends beyond just one simulation period. For example, if mission requirements mandate a 70$\%$ coverage requirement satisfaction, then our constellation meets the mandate.

Furthermore, Figure \ref{fig:14} shows that there exist several values of $\phi_0$ satisfying $>$ 90$\%$ of the coverage requirement. Figure \ref{fig:15} shows the coverage map for one such angle,  $\phi_0 = 162^{\circ}$. Not only is 97$\%$ of the coverage requirement satisfied, but the good ratio of green to black shows that robustness to delay is still present.

\begin{figure}[H]
  \centering
  \begin{minipage}[b]{0.45\textwidth}
    \includegraphics[width=\textwidth]{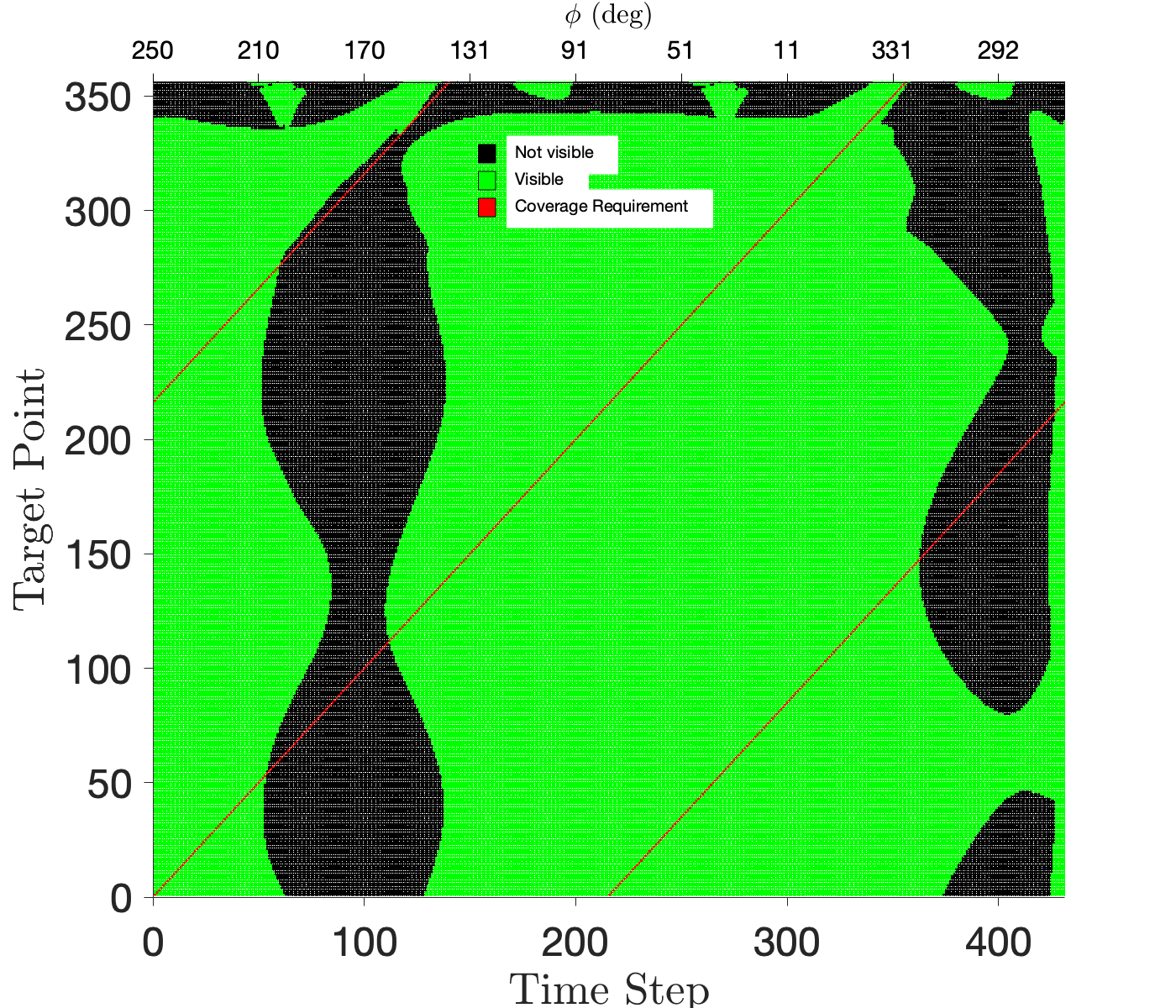}
    \caption{Coverage Robustness for 2 Departure Windows with $\phi_0 = 250^{\circ}$.  Constellation satisfies 71\% of coverage requirement. }
    \label{fig:16}
  \end{minipage}
  \hfill
  \begin{minipage}[b]{0.45\textwidth}
    \includegraphics[width=\textwidth]{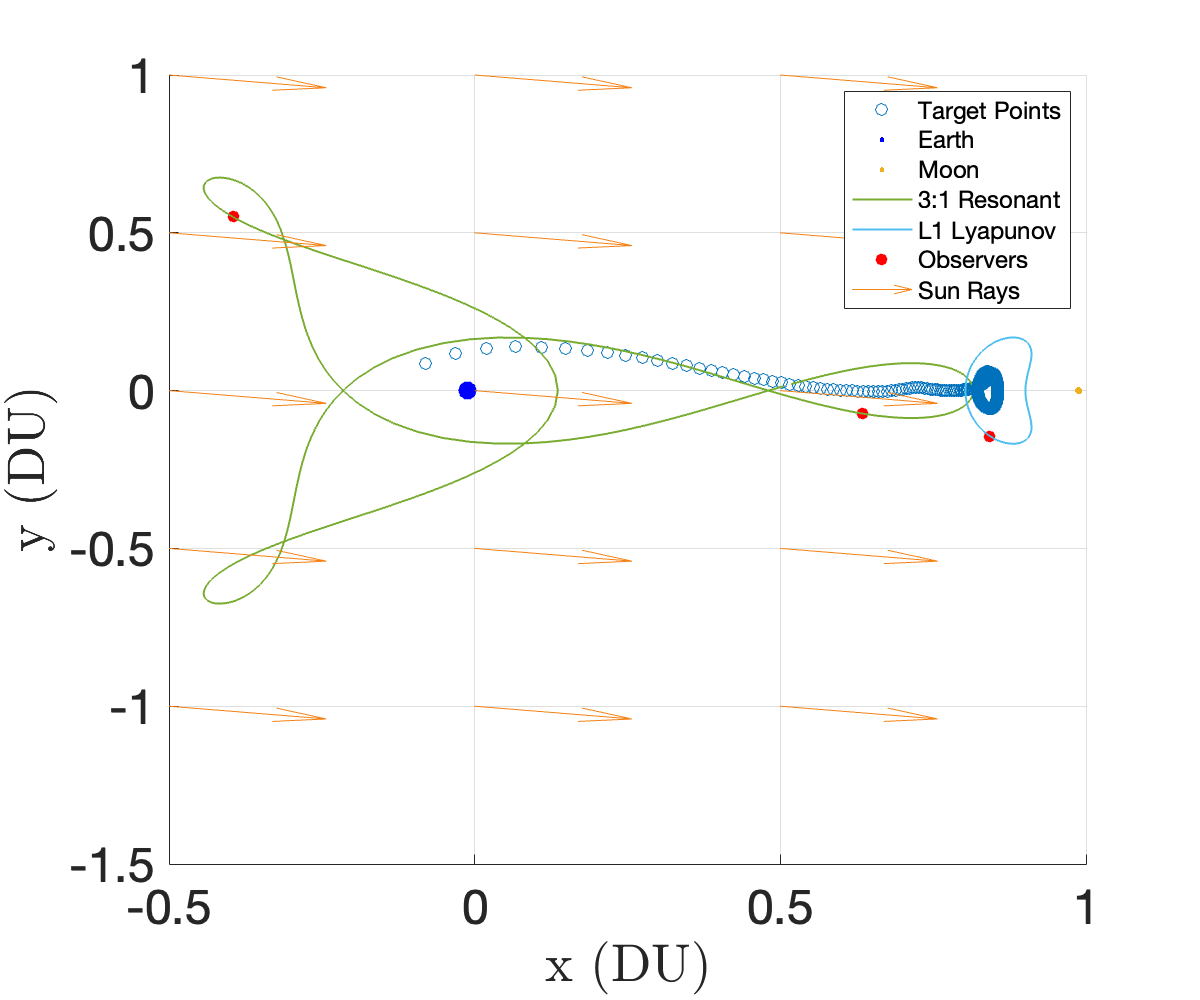}
    \caption{Satellite Positions at Time Step 100.  Sun rays are shown in orange. Most target points are backlit, resulting in poor illumination. }
    \label{fig:17}
  \end{minipage}
\end{figure}

 Of all the investigated $\phi_0$ in the range [0, 2$\pi$], $\phi_0 = 250^{\circ}$ was the worst performing constellation. However, it still satisfied 71 $\%$ of the coverage requirement. Figure \ref{fig:16} shows the coverage map for this choice of $\phi_0$, illustrating that despite only satisfying 71 $\%$ of the coverage requirement, the constellation still boasts good robustness to delays (shown via the ratio of green to black). There is a band of invisibility (appearing as an hourglass figure in black) in the approximate time step range of 70 to 130. This is due to unfavorable solar phase angles $\psi$ that are incurred in this region. Figure \ref{fig:17} shows the position of the observers as well as the direction of the Sun's rays at time step 100. We see that the positions of the Sun, target points, and observer on the Lyapunov orbit result in solar phase angles $>\frac{\pi}{2}$, and for some target points even closer to $\pi$. In this scenario the target is backlit and thus is poorly illuminated, making it difficult for optical sensors to monitor.

\section{Conclusions and Future Work}\label{sec5}
\subsection{Conclusion} 
 The main contribution of this work is extending the APC decomposition and ILP formulation by Lee et al. to analyze and optimize the cislunar SDA constellation design. We showed the developed formulation can optimize a constellation in cislunar space in response to spatio-temporally varying demand, such as tracking targets in an orbit with multiple departure windows (Case 1) or tracking multiple targets in different orbits simultaneously (Case 2). The robustness of the resulting solution against delays in departure and the initial sun phase angle was also analyzed.  We find that the optimal constellation size grows sublinearly with the number of departure windows. Also, we find that the constellation optimized for one simulation period also exhibits good performance past this time, despite changes in solar illumination. 
\subsection{Future Work} 
Multiple future work directions are possible. Instead of the Sun on a coplanar orbit, we may incline it by $\approx 5^{\circ}$ to reflect the Earth-Moon system's tilt with respect to the ecliptic. Furthermore, the Sun's gravitational perturbation can be incorporated to determine its effect on the stability and periodicity of candidate orbits in the CR3BP dynamics. Also, the objective of optimization can incorporate more complex costs. We can penalize constellations that place satellites on separate orbits and include a weighting term to reflect the degree of its importance. Lastly, the sensitivity of the solution architecture to the apparent magnitude cutoff is another direction for future work. Below a stringent enough cutoff (requiring targets to be unreasonably bright), the optimization program should become infeasible. Conversely, above a certain cutoff (assuming targets are visible when extremely dim) the optimization program should yield very small constellation sizes. 

\bmhead{Acknowledgment}
The authors gratefully acknowledge support for this research from the Air Force Office of Scientific Research (AFOSR), as part of the Space University Research Initiative (SURI), grant FA9550-22-1-0092 (grant principal investigator: J. Crassidis from University at Buffalo, The State University of New York).

\bmhead{Conflict of Interest Statement}
On behalf of all authors, the corresponding author states that there is no conflict of interest.

\bibliography{references}


\begin{thebibliography}{17}
\ifx \bisbn   \undefined \def \bisbn  #1{ISBN #1}\fi
\ifx \binits  \undefined \def \binits#1{#1}\fi
\ifx \bauthor  \undefined \def \bauthor#1{#1}\fi
\ifx \batitle  \undefined \def \batitle#1{#1}\fi
\ifx \bjtitle  \undefined \def \bjtitle#1{#1}\fi
\ifx \bvolume  \undefined \def \bvolume#1{\textbf{#1}}\fi
\ifx \byear  \undefined \def \byear#1{#1}\fi
\ifx \bissue  \undefined \def \bissue#1{#1}\fi
\ifx \bfpage  \undefined \def \bfpage#1{#1}\fi
\ifx \blpage  \undefined \def \blpage #1{#1}\fi
\ifx \burl  \undefined \def \burl#1{\textsf{#1}}\fi
\ifx \doiurl  \undefined \def \doiurl#1{\url{https://doi.org/#1}}\fi
\ifx \betal  \undefined \def \betal{\textit{et al.}}\fi
\ifx \binstitute  \undefined \def \binstitute#1{#1}\fi
\ifx \binstitutionaled  \undefined \def \binstitutionaled#1{#1}\fi
\ifx \bctitle  \undefined \def \bctitle#1{#1}\fi
\ifx \beditor  \undefined \def \beditor#1{#1}\fi
\ifx \bpublisher  \undefined \def \bpublisher#1{#1}\fi
\ifx \bbtitle  \undefined \def \bbtitle#1{#1}\fi
\ifx \bedition  \undefined \def \bedition#1{#1}\fi
\ifx \bseriesno  \undefined \def \bseriesno#1{#1}\fi
\ifx \blocation  \undefined \def \blocation#1{#1}\fi
\ifx \bsertitle  \undefined \def \bsertitle#1{#1}\fi
\ifx \bsnm \undefined \def \bsnm#1{#1}\fi
\ifx \bsuffix \undefined \def \bsuffix#1{#1}\fi
\ifx \bparticle \undefined \def \bparticle#1{#1}\fi
\ifx \barticle \undefined \def \barticle#1{#1}\fi
\bibcommenthead
\ifx \bconfdate \undefined \def \bconfdate #1{#1}\fi
\ifx \botherref \undefined \def \botherref #1{#1}\fi
\ifx \url \undefined \def \url#1{\textsf{#1}}\fi
\ifx \bchapter \undefined \def \bchapter#1{#1}\fi
\ifx \bbook \undefined \def \bbook#1{#1}\fi
\ifx \bcomment \undefined \def \bcomment#1{#1}\fi
\ifx \oauthor \undefined \def \oauthor#1{#1}\fi
\ifx \citeauthoryear \undefined \def \citeauthoryear#1{#1}\fi
\ifx \endbibitem  \undefined \def \endbibitem {}\fi
\ifx \bconflocation  \undefined \def \bconflocation#1{#1}\fi
\ifx \arxivurl  \undefined \def \arxivurl#1{\textsf{#1}}\fi
\csname PreBibitemsHook\endcsname

\bibitem[\protect\citeauthoryear{Folta et~al.}{2013}]{folta2013preliminary}
\begin{bchapter}
\bauthor{\bsnm{Folta}, \binits{D.C.}},
\bauthor{\bsnm{Pavlak}, \binits{T.A.}},
\bauthor{\bsnm{Haapala}, \binits{A.F.}},
\bauthor{\bsnm{Howell}, \binits{K.C.}}:
\bctitle{Preliminary design considerations for access and operations in earth-moon l1/l2 orbits}.
In: \bbtitle{AAS/AIAA Spaceflight Mechanics Meeting}
(\byear{2013})
\end{bchapter}
\endbibitem

\bibitem[\protect\citeauthoryear{Haapala et~al.}{2013}]{haapala2013trajectory}
\begin{bchapter}
\bauthor{\bsnm{Haapala}, \binits{A.}},
\bauthor{\bsnm{Vaquero}, \binits{M.}},
\bauthor{\bsnm{Pavlak}, \binits{T.}},
\bauthor{\bsnm{Howell}, \binits{K.C.}},
\bauthor{\bsnm{Folta}, \binits{D.C.}}:
\bctitle{Trajectory selection strategy for tours in the earth-moon system}.
In: \bbtitle{AAS/AIAA Astrodynamics Specialist Conference, Hilton Head, South Carolina}
(\byear{2013})
\end{bchapter}
\endbibitem

\bibitem[\protect\citeauthoryear{Cunio et~al.}{2020}]{cunio2020payload}
\begin{bchapter}
\bauthor{\bsnm{Cunio}, \binits{P.M.}},
\bauthor{\bsnm{Bever}, \binits{M.J.}},
\bauthor{\bsnm{Flewelling}, \binits{B.R.}}:
\bctitle{Payload and constellation design for a solar exclusion-avoiding cislunar ssa fleet}.
In: \bbtitle{Advanced Maui Optical and Space Surveillance Technologies Conference (AMOS)}
(\byear{2020})
\end{bchapter}
\endbibitem

\bibitem[\protect\citeauthoryear{Thompson et~al.}{2021}]{thompson2021cislunar}
\begin{bchapter}
\bauthor{\bsnm{Thompson}, \binits{M.R.}},
\bauthor{\bsnm{R{\'e}}, \binits{N.P.}},
\bauthor{\bsnm{Meek}, \binits{C.}},
\bauthor{\bsnm{Cheetham}, \binits{B.}}:
\bctitle{Cislunar orbit determination and tracking via simulated space-based measurements}.
In: \bbtitle{Advanced Maui Optical and Space Surveillance Technologies Conference, Maui, HI}
(\byear{2021})
\end{bchapter}
\endbibitem

\bibitem[\protect\citeauthoryear{Wilmer et~al.}{2022}]{wilmer}
\begin{barticle}
\bauthor{\bsnm{Wilmer}, \binits{A.P.}},
\bauthor{\bsnm{Bettinger}, \binits{R.A.}},
\bauthor{\bsnm{Little}, \binits{B.D.}}:
\batitle{Cislunar periodic orbits for earth–moon l1 and l2 lagrange point surveillance}.
\bjtitle{Journal of Spacecraft and Rockets}
\bvolume{59}(\bissue{6}),
\bfpage{1809}--\blpage{1820}
(\byear{2022})
\doiurl{10.2514/1.A35337}
{\href{https://arxiv.org/abs/https://doi.org/10.2514/1.A35337}{{https://doi.org/10.2514/1.A35337}}}
\end{barticle}
\endbibitem

\bibitem[\protect\citeauthoryear{Dahlke et~al.}{2022}]{dahlke2022preliminary}
\begin{bchapter}
\bauthor{\bsnm{Dahlke}, \binits{J.A.}},
\bauthor{\bsnm{Wilmer}, \binits{A.P.}},
\bauthor{\bsnm{Bettinger}, \binits{R.A.}}:
\bctitle{Preliminary comparative assessment of l2 and l3 surveillance using select cislunar periodic orbits}.
In: \bbtitle{AAS/AIAA Astrodynamics Specialist Conference},
pp. \bfpage{1}--\blpage{19}
(\byear{2022})
\end{bchapter}
\endbibitem

\bibitem[\protect\citeauthoryear{Zimovan et~al.}{2017}]{zimovan2017near}
\begin{bchapter}
\bauthor{\bsnm{Zimovan}, \binits{E.M.}},
\bauthor{\bsnm{Howell}, \binits{K.C.}},
\bauthor{\bsnm{Davis}, \binits{D.C.}}:
\bctitle{Near rectilinear halo orbits and their application in cis-lunar space}.
In: \bbtitle{3rd IAA Conference on Dynamics and Control of Space Systems, Moscow, Russia},
vol. \bseriesno{20},
p. \bfpage{40}
(\byear{2017})
\end{bchapter}
\endbibitem

\bibitem[\protect\citeauthoryear{Frueh et~al.}{2021}]{frueh2021cislunar}
\begin{bchapter}
\bauthor{\bsnm{Frueh}, \binits{C.}},
\bauthor{\bsnm{Howell}, \binits{K.}},
\bauthor{\bsnm{DeMars}, \binits{K.}},
\bauthor{\bsnm{Bhadauria}, \binits{S.}},
\bauthor{\bsnm{Gupta}, \binits{M.}}:
\bctitle{Cislunar space traffic management: Surveillance through earth-moon resonance orbits}.
In: \bbtitle{8th European Conference on Space Debris}
(\byear{2021})
\end{bchapter}
\endbibitem

\bibitem[\protect\citeauthoryear{Gupta et~al.}{2023}]{guptaconstellation}
\begin{bchapter}
\bauthor{\bsnm{Gupta}, \binits{M.}},
\bauthor{\bsnm{Howell}, \binits{K.C.}},
\bauthor{\bsnm{Frueh}, \binits{C.}}:
\bctitle{Constellation design to support cislunar surveillance leveraging sidereal resonant orbits}.
In: \bbtitle{33rd AAS/AIAA Space Flight Mechanics Meeting, Austin, Texas}
(\byear{2023})
\end{bchapter}
\endbibitem

\bibitem[\protect\citeauthoryear{Klonowski et~al.}{2023}]{Klonowski}
\begin{barticle}
\bauthor{\bsnm{Klonowski}, \binits{M.}},
\bauthor{\bsnm{Holzinger}, \binits{M.J.}},
\bauthor{\bsnm{Fahrner}, \binits{N.O.}}:
\batitle{Optimal cislunar architecture design using monte carlo tree search methods}.
\bjtitle{The Journal of the Astronautical Sciences}
\bvolume{70}(\bissue{3}),
\bfpage{17}
(\byear{2023})
\end{barticle}
\endbibitem

\bibitem[\protect\citeauthoryear{Vendl and Holzinger}{2021}]{vendl2021cislunar}
\begin{barticle}
\bauthor{\bsnm{Vendl}, \binits{J.K.}},
\bauthor{\bsnm{Holzinger}, \binits{M.J.}}:
\batitle{Cislunar periodic orbit analysis for persistent space object detection capability}.
\bjtitle{Journal of Spacecraft and Rockets}
\bvolume{58}(\bissue{4}),
\bfpage{1174}--\blpage{1185}
(\byear{2021})
\end{barticle}
\endbibitem

\bibitem[\protect\citeauthoryear{Visonneau et~al.}{2023}]{visonneau2023optimizing}
\begin{botherref}
\oauthor{\bsnm{Visonneau}, \binits{L.}},
\oauthor{\bsnm{Shimane}, \binits{Y.}},
\oauthor{\bsnm{Ho}, \binits{K.}}:
Optimizing multi-spacecraft cislunar space domain awareness systems via hidden-genes genetic algorithm.
The Journal of the Astronautical Sciences
\textbf{70}(22)
(2023)
\end{botherref}
\endbibitem

\bibitem[\protect\citeauthoryear{{Lee} et~al.}{2020}]{2020JSpRo..57.1309L}
\begin{barticle}
\bauthor{\bsnm{{Lee}}, \binits{H.W.}},
\bauthor{\bsnm{{Shimizu}}, \binits{S.}},
\bauthor{\bsnm{{Yoshikawa}}, \binits{S.}},
\bauthor{\bsnm{{Ho}}, \binits{K.}}:
\batitle{{Satellite Constellation Pattern Optimization for Complex Regional Coverage}}.
\bjtitle{Journal of Spacecraft and Rockets}
\bvolume{57}(\bissue{6}),
\bfpage{1309}--\blpage{1327}
(\byear{2020})
\doiurl{10.2514/1.A34657}
{\href{https://arxiv.org/abs/1910.00672}{{arXiv:1910.00672}}}
{[math.OC]}
\end{barticle}
\endbibitem

\bibitem[\protect\citeauthoryear{{Krag}}{1974}]{1974STIN...7512024K}
\begin{botherref}
\oauthor{\bsnm{{Krag}}, \binits{W.E.}}:
{Visible magnitude of typical satellites in synchronous orbits}
(1974)
\end{botherref}
\endbibitem

\bibitem[\protect\citeauthoryear{Vaquero and Senent}{2018}]{48975_2018}
\begin{botherref}
\oauthor{\bsnm{Vaquero}, \binits{M.}},
\oauthor{\bsnm{Senent}, \binits{J.}}:
{Poincare: A Multi-Body, Multi-System Trajectory Design Tool}
(2018)
\doiurl{2014/48975}
\end{botherref}
\endbibitem

\bibitem[\protect\citeauthoryear{Broucke}{1968}]{1968port.book.....B}
\begin{botherref}
\oauthor{\bsnm{Broucke}, \binits{R.A.}}:
Periodic orbits in the restricted three body problem with earth-moon masses.
Technical report
(1968)
\end{botherref}
\endbibitem

\bibitem[\protect\citeauthoryear{Rackauckas and Nie}{2017}]{rackauckas2017differentialequations}
\begin{botherref}
\oauthor{\bsnm{Rackauckas}, \binits{C.}},
\oauthor{\bsnm{Nie}, \binits{Q.}}:
Differential{E}quations.jl--a performant and feature-rich ecosystem for solving differential equations in {J}ulia.
Journal of Open Research Software
\textbf{5}(1)
(2017)
\end{botherref}
\endbibitem

\end{thebibliography}

\newpage
\begin{appendices}

\section{Optimizing the Constellation Pattern Vector}
Algorithm \ref{alg:1} presents pseudocode for the formulation. In this work, we utilize Gurobi MATLAB for computation. Computation is done with an M1 Pro CPU using 8 high-performance cores. Algorithm run time is $<$ 300 seconds for each $N \in\{1, 2, 4, 8, 16\}$.
\begin{algorithm}[H]
\caption{}\label{alg:1}
\begin{algorithmic}
\Require $\mathcal{Z}$, $\mathcal{J}$, $\boldsymbol{f}$, $g$, $M$, $T$, tstep

\State tspan $\gets$ 0:tstep:T
\State L $\gets$ {\fontfamily{qcr}\selectfont
length}(tspan)
\State $\mathcal{S} \gets$ {\fontfamily{qcr}\selectfont
suntrajectory}(tspan)
\\\hrulefill
\State V $\gets$ zeros(L, $\lvert \mathcal{Z} \rvert$, $\lvert \mathcal{J} \rvert$) \Comment{Precompute}
\ForEach{$z_j \in \mathcal{Z}$} 

\State $[\mathcal{R}, \mathcal{V}] \gets$ {\fontfamily{qcr}\selectfont
cr3bp}($z_j, \text{tspan}$) 

\ForEach{$(r_i, s_i) \in \{(r_m, s_n)\in \mathcal{R} \times \mathcal{S} \ \lvert \  n = m\}$}
\ForEach{$p_k \in \mathcal{J}$}

\State m $\gets g(p_k, r_i, s_i)$
\State V[i, j, k] $\gets$ m $\leq M\ $? 1 : 0

\EndFor
\EndFor

\EndFor
\\\hrulefill
\State $\bold{V}$ $\gets$ [ ] \Comment{Phase Orbits}
\For{$k = 1$ to $\lvert \mathcal{J} \rvert$}
\State C $\gets$ [ ]
\For{$j = 1$ to $\lvert \mathcal{Z} \rvert$}
\State A $\gets$ V[: , j , k]
\State B $\gets$ zeros(L, L)
\For{$i = 1$ to L}
\State B[: , $i$] $\gets$ {\fontfamily{qcr}\selectfont
circshift}(A, $i-1$)
\EndFor
\State C $\gets$ {\fontfamily{qcr}\selectfont
horzcat}(C, B)
\EndFor
\State $\bold{V}$ $\gets$ {\fontfamily{qcr}\selectfont
vertcat}($\bold{V}$, C)
\EndFor

\\\hrulefill

\State c $\gets$  ones($\lvert \mathcal{Z} \rvert$ * L) \Comment{Optimize}

\State result $\gets$ {\fontfamily{qcr}\selectfont
integer-linear-program}($\bold{V}$, $\boldsymbol{f}$, c)
\end{algorithmic}
\end{algorithm}
\newpage
\section{CR3BP System Parameters and Candidate Orbit Initial Conditions}
System parameters for the Earth-Moon CR3BP with the Sun as illumination source are presented in Table \ref{tab:1} and Table \ref{tab:2}. Initial conditions for candidate orbits are shown below in Table \ref{tab:3}. Orbits were generated via a single shooting differential corrections algorithm written in Julia. We use the Vern7 ODE propagator (available in the Julia DifferentialEquations \cite{rackauckas2017differentialequations} module) in this algorithm. 

\begin{table}[h]
\caption{CR3BP Parameters}
\begin{tabular}{ll}
\multicolumn{2}{c}{}\\ \hline
Parameter                        & Value        \\ \hline
Mass parameter, $\mu$ & 1.215058560962404e-02     \\
Length scale, DU (km)            & 384400.0     \\
Time scale, TU (s)               & 3.751902619517228e+05 \\
Relative distance of Sun (DU)    & 389.17794    \\
Angular velocity of Sun (rad/TU) & -0.9253018261815922 \\ 
Earth Radius, (km) & 6371 \\ 
Moon Radius, (km) & 1737.4
\end{tabular}
\label{tab:1}
\end{table}

\begin{table}[h]
\caption{Illumination Parameters}
\begin{tabular}{ll}
\multicolumn{2}{c}{}   \\ \hline
Parameter                        & Value        \\ \hline
Sun apparent magnitude, $m_{\text{sun}}$  & -26.74     \\
Target specular reflection coef., $a_{\text{spec}}$         & 0.0     \\
Target diffuse reflection coef., $a_{\text{diff}}$            & 0.2 \\
Target Diameter, $d$  (km)    & 0.001 
\end{tabular}
\label{tab:2}
\end{table}

\begin{table}[h]
\caption{Candidate Orbit Initial Conditions}
\begin{tabular}{llll}
\multicolumn{4}{c}{}   \\ \hline
& 3:1 Resonant  &  2:1 Resonant & L1 Lyapunov  \\ \hline
$x_0$ &0.13603399956670137  &    0.9519486347314083 & 0.65457084231188\\ 
$y_0$ & 0 & 0 & 0\\
$z_0$ & 0 & 0 & 0\\
$\dot{x}_0$ &1.9130717669166003e-12 & 0 & 3.887957091335523e-13\\
$\dot{y}_0$ &3.202418276067991 & -0.952445273435512 & 0.7413347560791179\\
$\dot{z}_0$ &0 & 0 & 0\\ 
T & 6.45 & 6.45 & 6.45 \\ \hline
&L2 Lyapunov & L1 Lyapunov (short) & L2 Halo (short) \\ \hline
$x_0$ & 0.9982702689023665 & 0.8027692908754149 & 1.1540242813087864\\ 
$y_0$  & 0 & 0 & 0\\ 
$z_0$  & 0 & 0 & -0.1384196144071876\\
$\dot{x}_0$ & -2.5322340091977996e-14 & -1.1309830924549648e-14 & 4.06530060663289e-15 \\
$\dot{y}_0$ & 1.5325475708886613 & 0.33765564334938736 & -0.21493019200956867\\ 
$\dot{z}_0$ & 0 & 0 & 8.48098638414804e-15\\ 
T & 6.45 & 3.225 & 3.225
\end{tabular}
\label{tab:3}
\end{table}

\newpage
\section{Coverage Requirement Generation}
Algorithm \ref{alg:2} provides pseudocode for the coverage requirement satisfying the criteria outlined in section \ref{subsubsec4}. 
\begin{algorithm}[H]
\caption{}\label{alg:2}
\begin{algorithmic}
\Require $N$, $T$ 
\Comment{Number of Departure Windows (N) and Timesteps (T)}

\State $\boldsymbol{f}_1$ $\gets$ {\fontfamily{qcr}\selectfont
zeros}(T)

\\

\State step $\gets$ T
\State power $\gets$ 0
\State indices $\gets$ [1]

\While {power $<$ {\fontfamily{qcr}\selectfont
log$_2$}(N)}

\State step $\gets$ {\fontfamily{qcr}\selectfont floor}(step/2)
\State power $\gets$ power + 1
\State newindices $\gets$ {\fontfamily{qcr}\selectfont
zeros}(1, {\fontfamily{qcr}\selectfont
length}(indices))

\For {i = 1 to {\fontfamily{qcr}\selectfont
length}(indices)}
\State newindices[i] $\gets$ indices[i] + step
\EndFor

\State indices $\gets$ [indices, newindices]

\EndWhile
\\

\State $\boldsymbol{f}_1$[indices] $\gets$ 1

\Return $\boldsymbol{f}_1$

\end{algorithmic}
\end{algorithm}

\newpage
\section{Extended Results}
\begin{figure}[h]
    \centering
    \begin{tabular}{c|c}
        \includegraphics[width=.45\linewidth,valign=m]{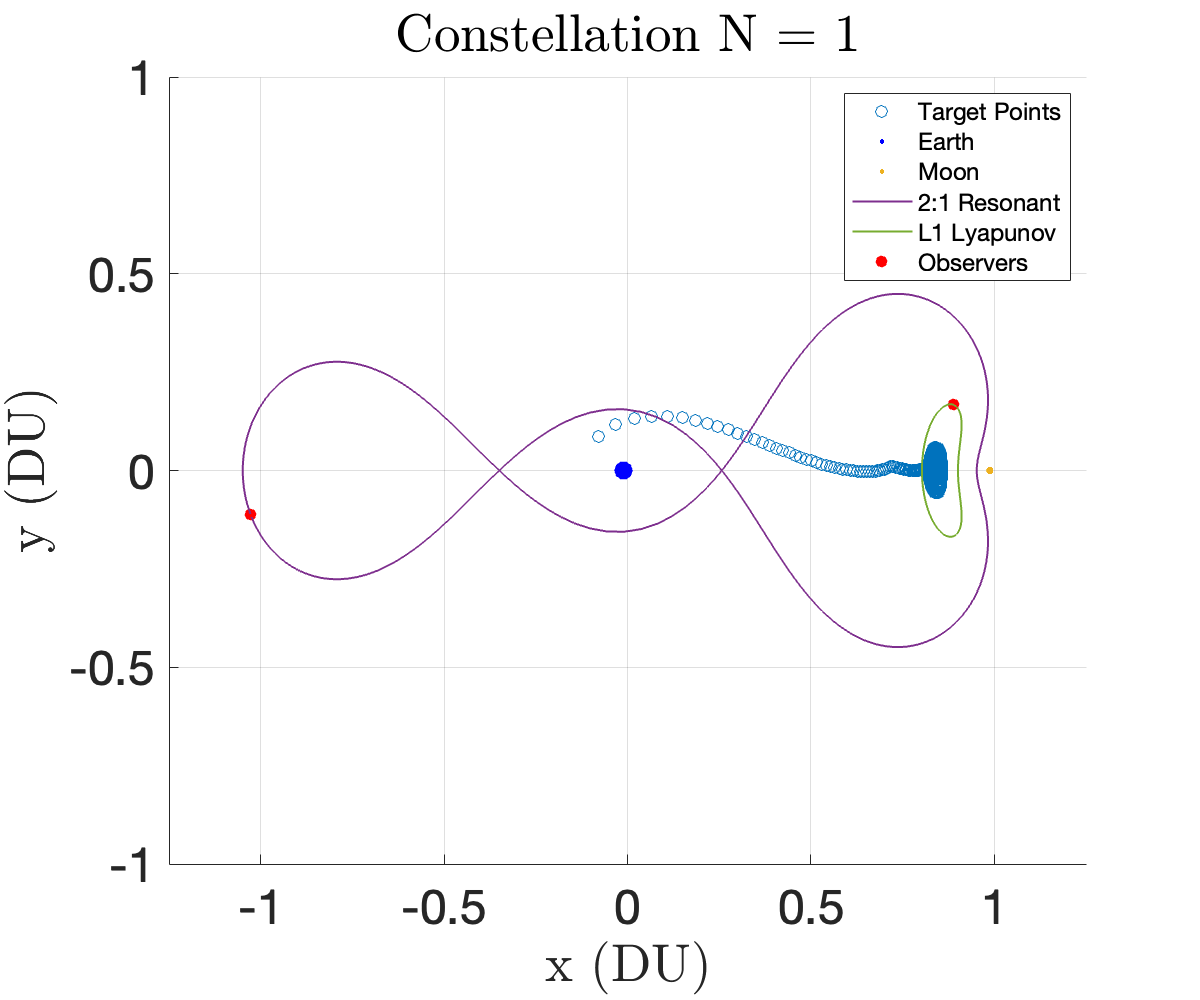}  & \includegraphics[width=.45\linewidth,valign=m] {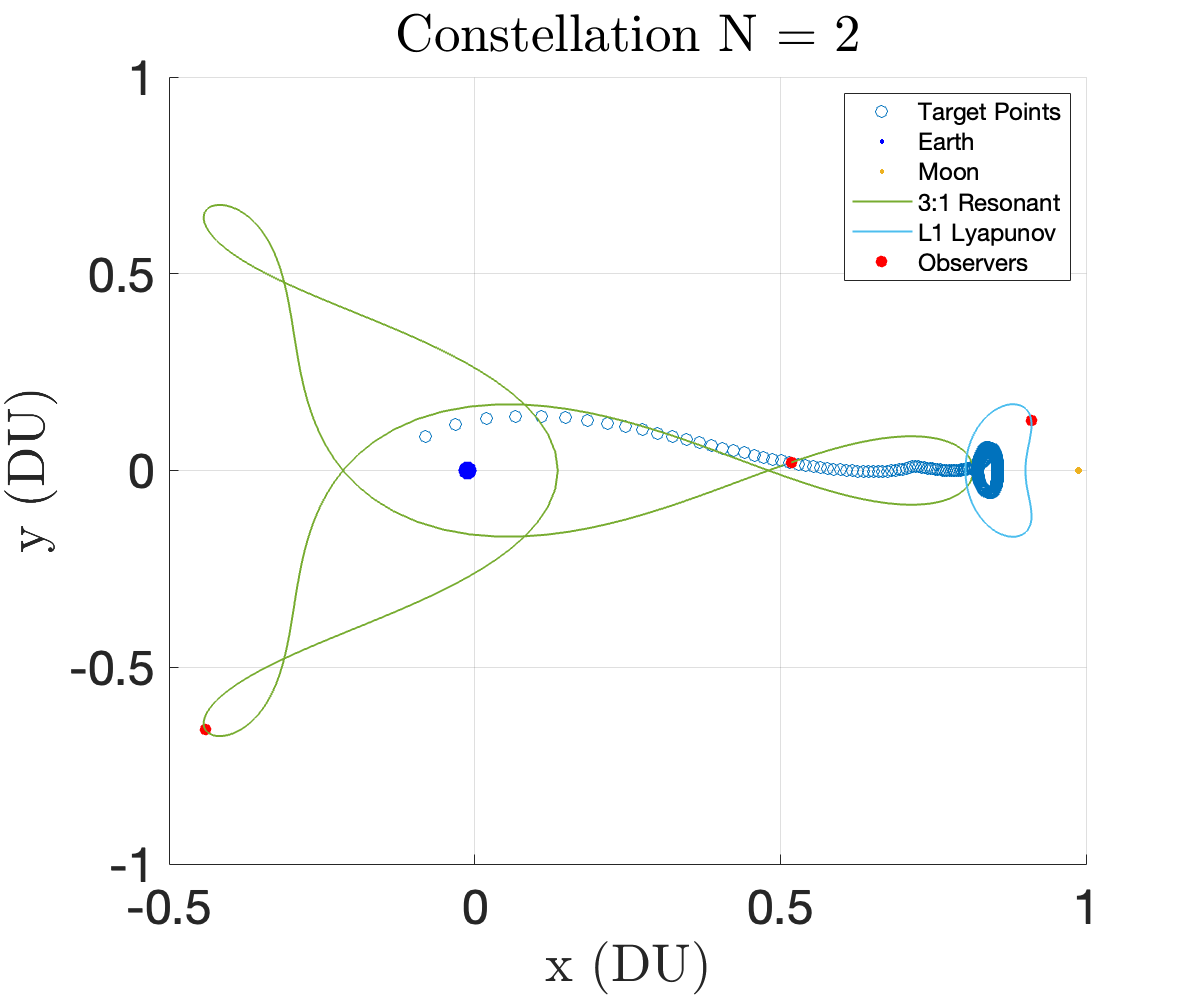} \\
        \\
        \hline
        \includegraphics[width=.45\linewidth,valign=m]{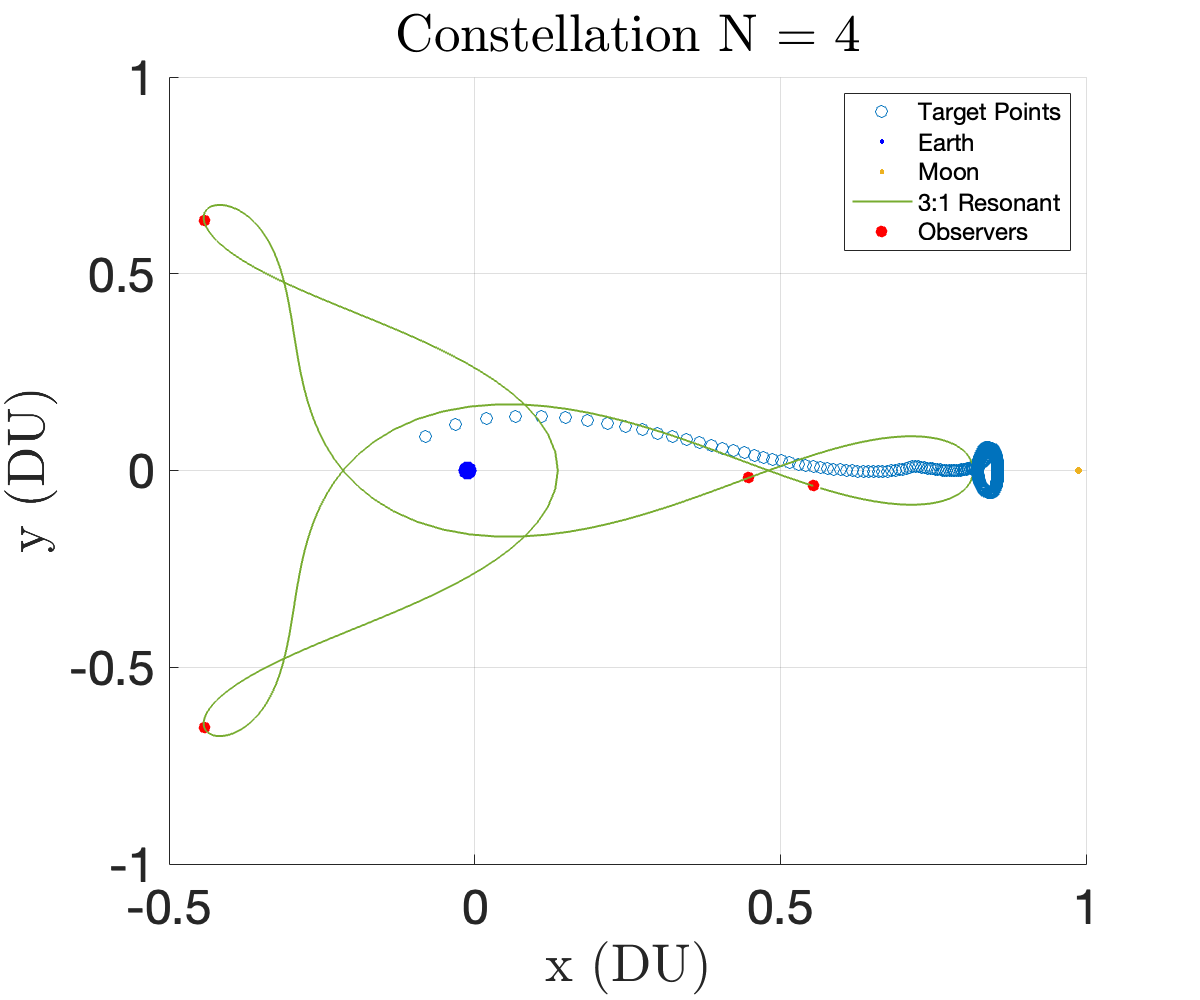} & \includegraphics[width=.45\linewidth,valign=m] {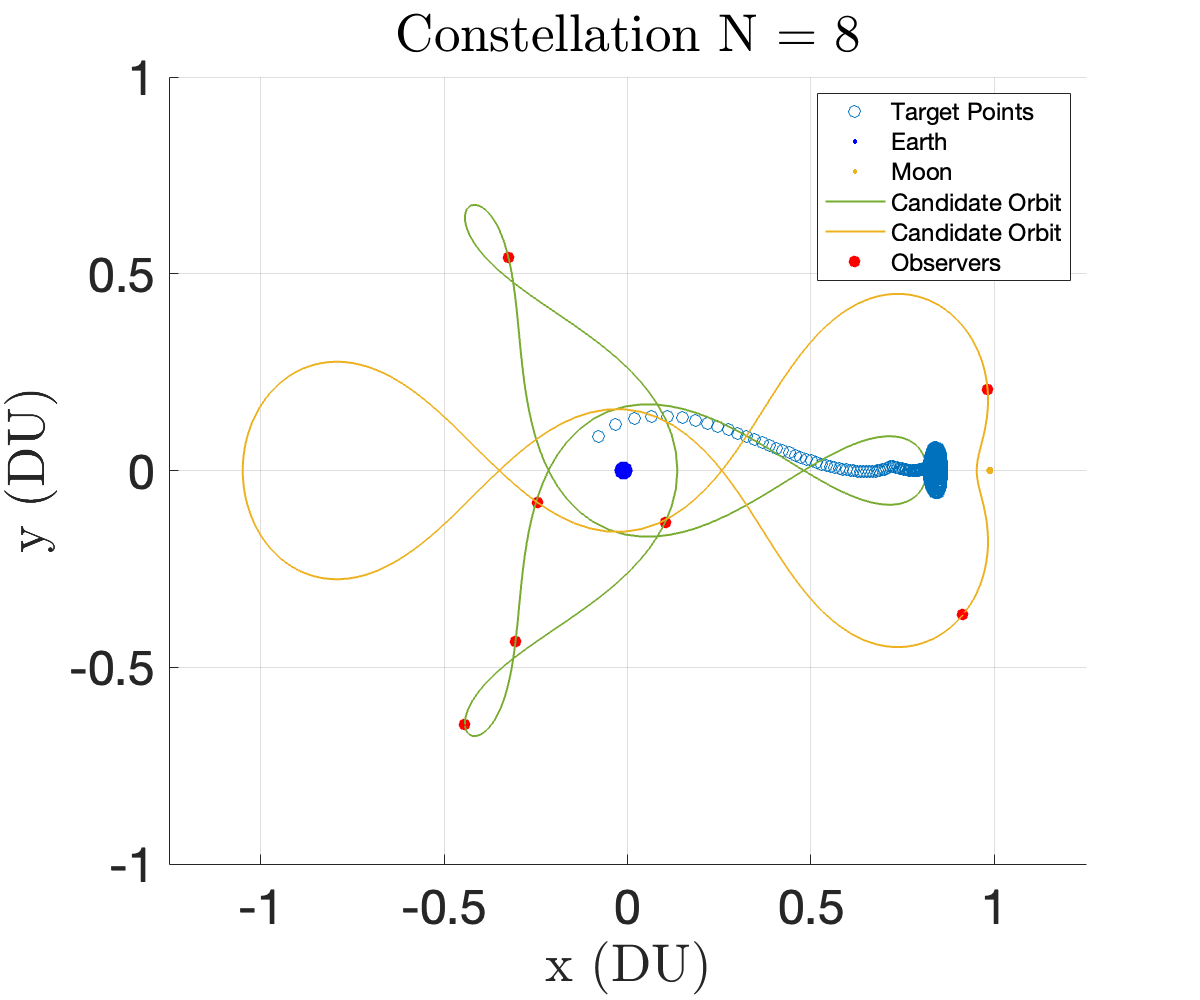} \\
        \\
        \hline
        \includegraphics[width=.45\linewidth,valign=m]{Figures/N16.png} & \includegraphics[width=.45\linewidth,valign=m]{Figures/windows_vs_sat.png}
    \end{tabular}
    \caption{Optimized Constellations for Given Number of Departure Windows}
    \label{fig:18}
\end{figure}

\begin{figure}[h]
    \centering
    \begin{tabular}{c|c}
        \includegraphics[width=.45\linewidth,valign=m]{Figures/robustnessN1.png}  & \includegraphics[width=.45\linewidth,valign=m] {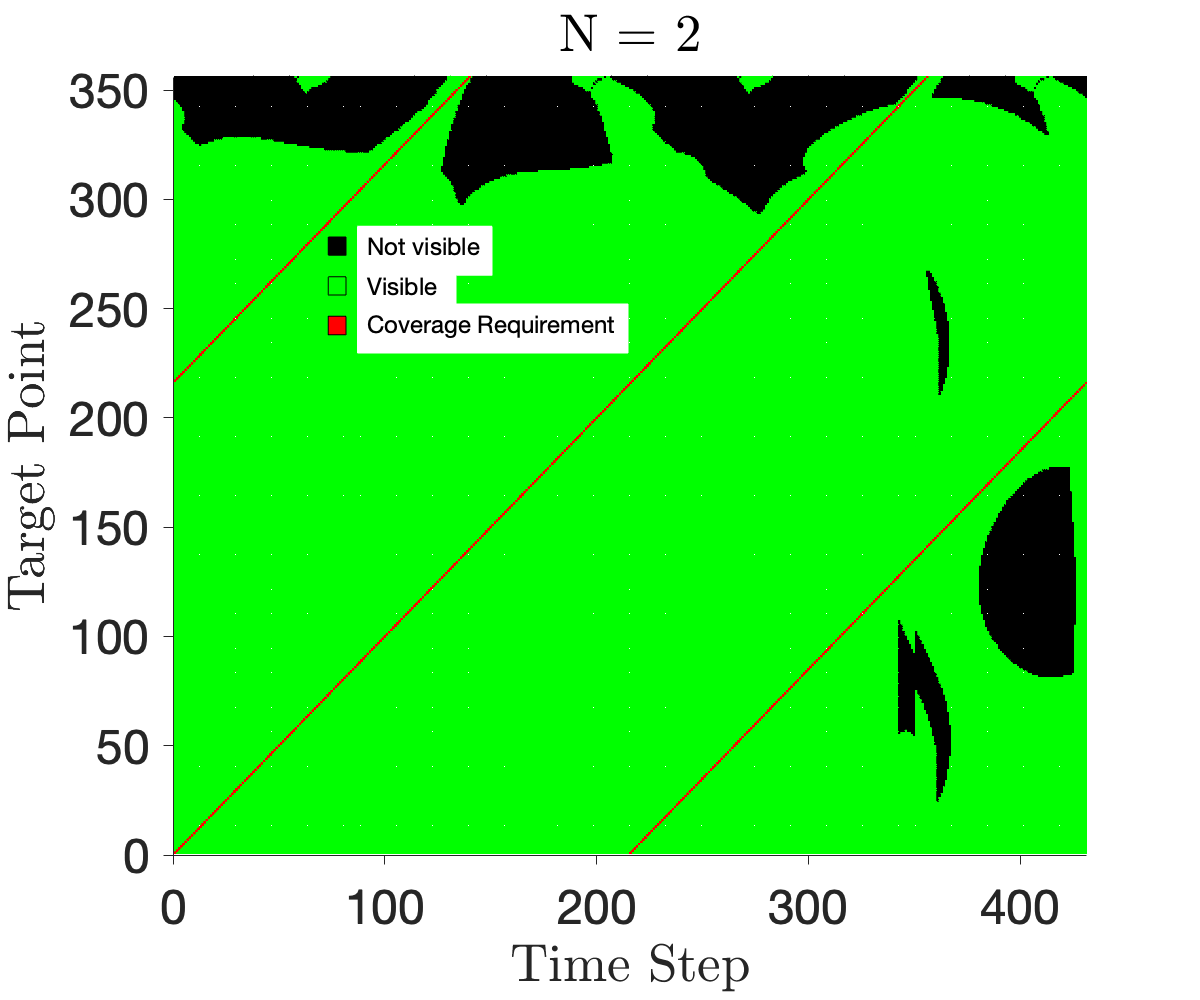} \\
        \\
        \hline
        \\
        \includegraphics[width=.45\linewidth,valign=m]{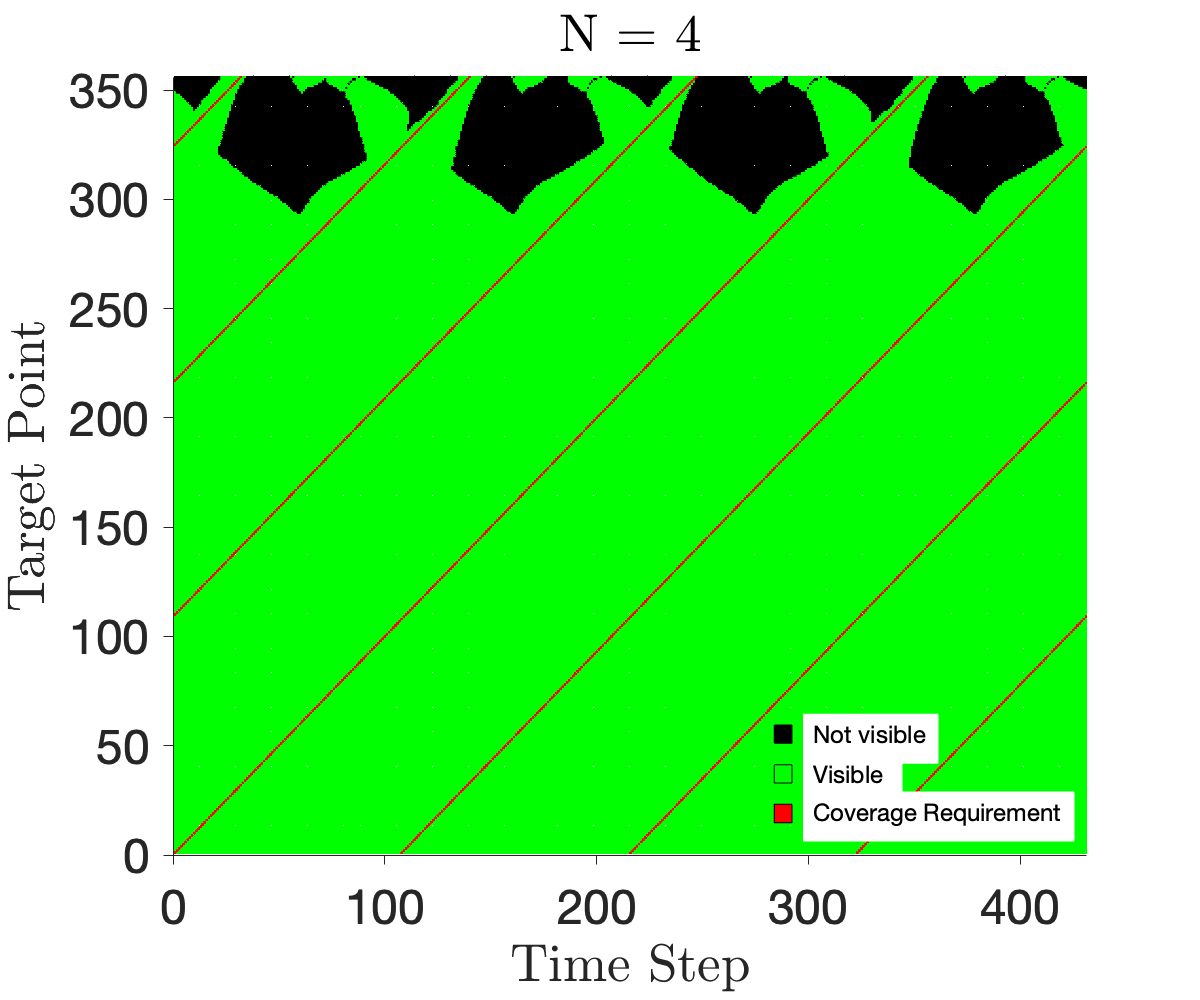} & \includegraphics[width=.45\linewidth,valign=m] {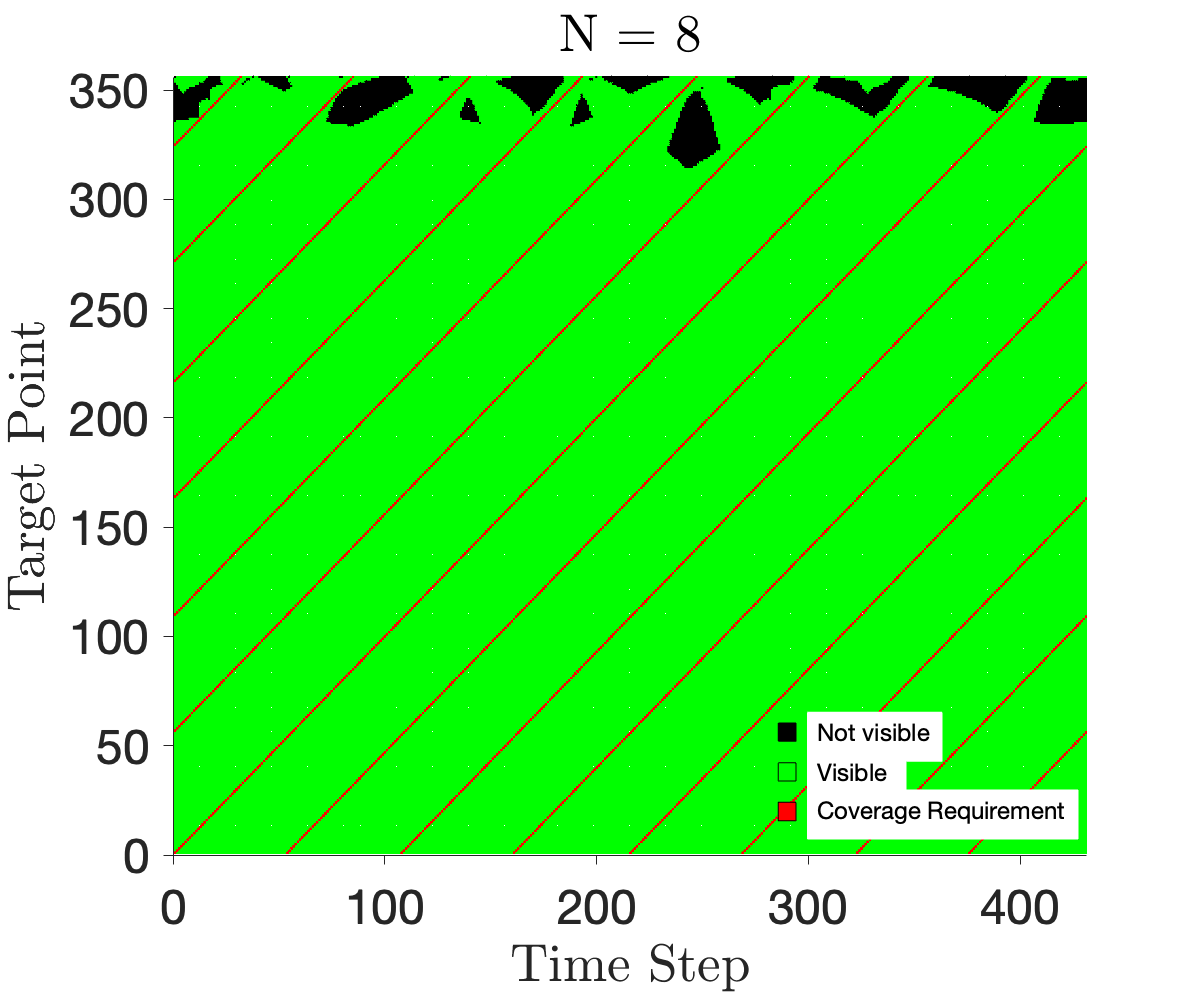} \\
        \\
        \hline
        \\
        \includegraphics[width=.45\linewidth,valign=m]{Figures/robustnessN16.png}

    \end{tabular}
    \caption{Constellation Robustness Binary Maps}
    \label{fig:19}
\end{figure}

Figures \ref{fig:18} and \ref{fig:19} show the constellations and indicator maps of all the optimization runs for each $N \in \{1, 2, 4, 8, 16\}$. Table \ref{tab:4} shows the phasing of the constellations relative to seed satellites.

\begin{sidewaystable}[h]
    \caption{Constellation Phasing Relative to Seed Satellites. Phasing is defined as a circular shift operation on seed satellite's position state history.}
    \centering
    \begin{tabular}{l|c|c|c|c|c}
    \multicolumn{6}{c}{} \\ 
    & N = 1  &  N = 2 & N = 4 & N = 8 & N = 16 \\ \hline 
    & & & & &\\ 
    3:1 Resonant &  & [62, 271] & [60, 164, 276, 377] & [4, 58, 114, 327] & [134, 162, 350, 380, 403] \\ 
    & & & & &\\ 
    2:1 Resonant & [230] & & & [21, 125, 390] & [71, 216, 256, 389] \\ 
    & & & & &\\ 
    L1 Lyapunov  & & & & & \\ 
    & & & & &\\ 
    L2 Lyapunov & & & & & \\ 
    & & & & &\\ 
    L1 Lyapunov (short) & [372] & [350] &  & & \\ 
    & & & & &\\ 
    L2 Halo (short) & & & & & \\ 
    
    \end{tabular}
    \label{tab:4}
\end{sidewaystable}

\newpage

\end{appendices}


\end{document}